\documentclass[12pt,a4paper]{amsart}
\usepackage{amssymb,enumerate,comment,amsmath}
\usepackage[table]{xcolor}

\usepackage[pagebackref]{hyperref}

\newcommand{\FF}{{\mathbb{F}}}
\newcommand{\QQ}{{\mathbb{Q}}}
\newcommand{\ZZ}{{\mathbb{Z}}}

\newcommand{\bC} {\mathbf C}
\newcommand{\bE} {\mathbf E}
\newcommand{\bF} {\mathbf F}
\newcommand{\bN} {\mathbf N}

\newcommand{\bG}{\mathbf G}
\newcommand{\bO}{{\mathbf O}}
\newcommand{\bZ}{\mathbf Z}

\newcommand{\cG} {\mathcal{G}}
\newcommand{\cP} {\mathcal{P}}
\newcommand{\cJ} {\mathcal{J}}

\newcommand{\fA} {\mathfrak A}
\newcommand{\fS} {\mathfrak S}

\newcommand{\alt}{\fA}
\newcommand{\sym}{\fS}

\newcommand{\ab}{{\operatorname{ab}}}
\newcommand{\ad}{{\operatorname{ad}}}
\newcommand{\Aut}{{\operatorname{Aut}}}
\newcommand{\cd}{{\operatorname{cd}}}
\newcommand{\diag}{{\operatorname{diag}}}
\newcommand{\Gal}{{\operatorname{Gal}}}
\newcommand{\Irr}{{\operatorname{Irr}}}
\newcommand{\Out}{{\operatorname{Out}}}
\newcommand{\res}{{\operatorname{res}}}
\newcommand{\St}{{\operatorname{St}}}
\newcommand{\Syl}{{\operatorname{Syl}}}
\newcommand{\type}{\operatorname}

\newcommand{\GL}{\operatorname{GL}}
\newcommand{\SL}{\operatorname{SL}}
\newcommand{\PGL}{\operatorname{PGL}}
\newcommand{\PSL}{\operatorname{PSL}}
\newcommand{\GU}{\operatorname{GU}}
\newcommand{\SU}{\operatorname{SU}}

\newcommand{\PSU}{\operatorname{PSU}}
\newcommand{\Spin}{\operatorname{Spin}}
\newcommand{\PSO}{\operatorname{PSO}}

\def\zent#1{{\bf Z}(#1)}
\newcommand{\cent} {\mathbf C}
\newcommand{\norm} {\mathbf N}
\newcommand{\wt}{\widetilde}
\newcommand{\tw}[1]{{}^#1\!}
\newcommand{{\tchi}}{\tilde\chi}
\newcommand{{\wG}}{{\widehat{G}}}
\newcommand{{\wpsi}}{{\widehat{\psi}}}
\renewcommand{\pmod}[1]{~({\rm mod}~#1)}

\let\eps=\epsilon
\let\veps=\varepsilon
\let\ga=\gamma
\let\la=\lambda
\let\Om=\Omega
\let\vhi=\varphi
\let\ze=\zeta

\newtheorem{thm}{Theorem}[section]
\newtheorem{lem}[thm]{Lemma}
\newtheorem{cor}[thm]{Corollary}
\newtheorem{prop}[thm]{Proposition}

\newtheorem{thmA}{Theorem}
\newtheorem{conjA}[thmA]{Conjecture}

\theoremstyle{definition}
\newtheorem{rem}[thm]{Remark}

\newtheorem{exmp}[thm]{Example}

\numberwithin{equation}{section}

\begin{document}

%%%%%%%%%%%%%%%%%%%%%%%%%%%%%%%%%%%%%%%%%%%%%%%%%%%%%%%%%%%%%
\title{A Brauer--Galois height zero conjecture}
%%%%%%%%%%%%%%%%%%%%%%%%%%%%%%%%%%%%%%%%%%%%%%%%%%%%%%%%%%%%%

\author{Gunter Malle}
\address[G. Malle]{FB Mathematik, RPTU, Postfach 3049,
  67653 Kaisers\-lautern, Germany.}
\email{malle@mathematik.uni-kl.de}

\author{Alexander Moret\'o}
\address[A. Moret\'o]{Departament de Matem\`atiques, Universitat de Val\`encia,
 46100 Burjassot, Val\`encia, Spain}
 \email{alexander.moreto@uv.es}

\author{Noelia Rizo}
\address[N. Rizo]{Departament de Matem\`atiques, Universitat de Val\`encia,
 46100 Burjassot, Val\`encia, Spain}
 \email{Noelia.Rizo@uv.es}

\author{A. A. Schaeffer Fry}
\address[A. A. Schaeffer Fry]{Dept. Mathematics - University of Denver, Denver
 CO 80210, USA; and
 Dept. Mathematics and Statistics - MSU Denver, Denver, CO 80217, USA}
\email{mandi.schaefferfry@du.edu}

\begin{abstract}
Recently, Malle and Navarro obtained a Galois strengthening of Brauer's height
zero conjecture for principal $p$-blocks when $p=2$, considering a particular
Galois automorphism of order~$2$. In this paper, for any prime $p$ we consider
a certain elementary abelian $p$-subgroup of the absolute Galois group and propose a
Galois version of Brauer's height zero conjecture for principal $p$-blocks.
We prove it when $p=2$ and also for arbitrary $p$ when $G$ does not involve
certain groups of Lie type of small rank as composition factors. Furthermore,
we prove it for almost simple groups and for $p$-solvable groups.
\end{abstract}

\thanks{
We thank Gabriel Navarro for suggesting Theorem \ref{thm:evenorder2},
for proposing to replace $\Om$ by $\cJ$ and for further helpful comments.
The first author gratefully acknowledges financial support by the DFG
-- Project-ID 286237555.
The second and third authors are supported by Ministerio de Ciencia e
Innovaci\'on (Grant PID2019-103854GB-I00 and Grant PID2022-137612NBI00 funded
by MCIN/AEI/10.13039/501100011033 and “ERDF A way of making Europe”) and a
CDEIGENT grant CIDEIG/2022/29 funded by Generalitat Valenciana. The second
author also acknowledges support by Generalitat Valenciana CIAICO/2021/163.
The fourth author gratefully acknowledges support from the National Science
Foundation, Award No. DMS-2100912, and her former institution, Metropolitan
State University of Denver, which holds the award and allows her to serve as PI} 



\date{\today}

\maketitle

%\pagestyle{myheadings}
%\markboth{for personal use only}{preliminary}

%%%%%%%%%%%%%%%%%%%%%%%%%%%%%%%%%%%%%%%%%%%%%%%%%%%%%%%%%%%%%%%%%%%%%%%%%
\section{Introduction}

Some of the fundamental local/global counting conjectures in representation
theory of finite groups were strengthened by Navarro \cite{N04} using the
action of Frobenius elements of $\Gal(\QQ^\ab/\QQ)$. This strengthening
has had a large impact on research in the last two decades. It was asked in
\cite{MN} whether there is also a strengthening of Brauer's height zero
conjecture (which was very recently solved \cite{MNST}) to include Galois
automorphisms. Brauer's conjecture, posed in 1955, asserts that if $G$ is a
finite group, and $B$ is a Brauer $p$-block of $G$ with defect group $D$, then
$D$ is abelian if and only if all the complex irreducible characters in $B$
have height zero. (Recall that $\chi\in\Irr(B)$ has height zero if
$\chi(1)_p=|G|_p/|D|$ is minimal possible.)  The main result of \cite{MN}
provides such a strengthening for the principal $2$-block. It was proved there
that if $\sigma$ denotes the Galois automorphism that fixes $2$-power order
roots of unity and complex-conjugates odd-order roots of unity, then all
$\sigma$-invariant irreducible characters in the principal $2$-block $B_0(G)$
of a finite group~$G$ have
$2'$-degree if and only if $G$ has an abelian Sylow $2$-subgroup. As pointed
out in \cite{MN}, this result does not extend to arbitrary $2$-blocks.

For $p$ a fixed prime number let $\cJ$ be the subgroup of $\Gal(\QQ^\ab/\QQ)$
consisting of the automorphisms of order~$p$ that fix all $p$-power order roots
of unity. We write $\Irr_\cJ(G)$ to denote the set of $\cJ$-invariant complex
irreducible characters of a finite group $G$ and $\Irr_\cJ(B_0(G))$ for
the $\cJ$-invariant irreducible characters in the principal $p$-block~$B_0(G)$. 

We propose the following conjecture.

\begin{conjA}   \label{galcon}
 Let $G$ be a finite group, let $p$ be a prime number and let $P\in\Syl_p(G)$.
 Then all characters in $\Irr_\cJ(B_0(G))$ have height zero if and only if
 $P$ is abelian.
\end{conjA}

This provides the desired Galois strengthening of the height zero conjecture
for principal $p$-blocks. We prove this conjecture when $p=2$ and also when $p$
is odd and certain simple groups are not composition factors of $G$.  Since
$\sigma\in\cJ$, this result also strengthens the main result of \cite{MN} when
$p=2$. (We remark that the original height zero conjecture had remained open
even for principal $p$-blocks until recently \cite{MN21}.)

\begin{comment}
Let $\cG=\cG(n)=\Gal(\QQ_n/\QQ)$. Then
$$\cG=\prod_{q}\cG_q,$$
where $$\cG_q=\Gal(\QQ_n/\QQ_{n_{q'}})\cong\Gal(\QQ_{n_q}/\QQ),$$ and the
(direct) product runs over the primes $q$ that divide $n$. Recall that 
$\cG_q$ is a cyclic group of order $\vhi(n_q)$ when $q>2$ or $q=2$ and
$n_2\leq4$ and $\cG_2\cong C_{n_2/4}\times C_2$ when $n_2\geq 8$.
Let $\cP_q\in\Syl_p(\cG_q)$. Note that $\cP_q$ is cyclic if $p>2$ or $q>2$ and
has rank~$2$ if $p=q=2$ and $n_2>4$. Let $\cP\in\Syl_p(\cG)$. We define
$$\Om=\Om(n)=\Om_1(\prod_{q\neq p}{\cP_q}),$$
where given a $p$-group $P$, $\Om_1(P)$ is the subgroup generated by the elements of order $p$ of~$P$.
\end{comment}

As pointed out in \cite{MN}, the ``only if'' implication of Conjecture
\ref{galcon} fails for non-principal blocks when $p=2$. It is also easy to find
examples in GAP \cite{gap} showing that this fails for non-principal $p$-blocks
also when $p$ is odd. On the other hand, the conclusion
of Theorem~\ref{main} is definitely false for $p=2$ if we define $\cJ$ to be
the set of elements of order $p$ of $\Gal(\QQ^\ab/\QQ)$. For instance, the
non-linear characters of the non-abelian 2-group ${\tt SmallGroup}(16, 6)$ in
\cite{gap} are not fixed by complex conjugation.  As we will see in
Remark~\ref{alt}, there are also counterexamples for odd primes. It is possible
that Conjecture~\ref{galcon} is true if we replace $\cJ$ by the set of elements
of $\Gal(\QQ^\ab/\QQ)$ that fix $p$-power roots of unity and have order a power
of~$p$. However, more non-abelian simple groups would appear as counterexamples
in the corresponding version of Theorem~\ref{im} (see Example~\ref{example:J1}).
 
The following is our main result.

\begin{thmA}   \label{main}
 Let $G$ be a finite group, let $p$ be a prime number and let $P\in\Syl_p(G)$.
 Suppose that $G$ does not have composition factors isomorphic to $S$ with
 $(S,p)$ a pair as listed in Theorem \ref{thm:as}. 
 Then Conjecture~\ref{galcon} holds for $G$. 
\end{thmA}

Our proof of Theorem \ref{main} relies on a partial strengthening of another
celebrated theorem in character theory of finite groups: the It\^o--Michler
theorem (see e.g.~\cite[Thm~7.1]{N18}).

\begin{thmA}   \label{im}
 Let $G$ be a finite group and let $p$ be a prime. Suppose that $G$ does not
 have composition factors isomorphic to $S$ with $(S,p)$ a pair as listed in
 Theorem~\ref{thm:as}. Then all characters in $\Irr_\cJ(G)$
 have $p'$-degree if and only if $G$ has a normal abelian Sylow $p$-subgroup.
 In particular, if $G$ is $p$-solvable then all characters in $\Irr_\cJ(G)$
 have $p'$-degree if and only if $G$ has a normal abelian Sylow $p$-subgroup.
\end{thmA}

The conclusion of Theorem \ref{im} is definitely false when $G=S$ for some
of the pairs $(S,p)$ listed in Theorem \ref{thm:as}, see Example~\ref{ex:PSL}.
Notice that in the $p$-solvable case we have a complete Galois version of the
It\^o--Michler theorem. This could be compared with Theorem~A of \cite{gri}.

We remark that we had reduced the proof of Theorem \ref{main} to two sets of
questions on simple groups, addressed in Section~\ref{sec:simple}, when the
preprint \cite{Gr}
appeared.  In that paper, a different extension of the It\^o--Michler theorem
is proved. Those of our questions relating to Theorem~\ref{im} are now
handled by a careful analysis and extension of the result of Theorem~C from
\cite{Gr}, while those relevant for the proof of Theorem~\ref{main}, involving
statements related to principal blocks, have to be solved differently.  As one
step, in particular, we show Conjecture~\ref{galcon} for almost simple groups.
We then prove Theorem~\ref{im} and the $p$-solvable case of
Conjecture~\ref{galcon} in Section~\ref{sec:IM} and use it in
Section~\ref{sec:bhz} to prove Theorem~\ref{main}.

%%%%%%%%%%%%%%%%%%%%%%%%%%%%%%%%%%%%%%%%%%%%%%%%%%%%%%%%%%%%%%%%%%%%%%%%%
\section{Almost simple and quasi-simple groups}   \label{sec:simple}
In this section we collect several results on characters of simple and
almost simple groups needed for the proofs of our main theorems.

Let $n$ be  an integer. We write $\QQ_n:=\QQ(\ze)$ where $\ze$ is a primitive
$n$th root of unity. In the proofs we will work with the following subgroup of
$\Gal(\QQ_n/\QQ)$.  Let $\Om=\Om(n)$ be the set of elements of order $p$ of
$\Gal(\QQ_n/\QQ_{n_p})$. Note that this definition depends on the integer $n$
and of course on $p$. If $G$ is a finite group, we often write $\Om$ for
$\Om(|G|)$. We let $\Irr_\Om(G)$
denote the set of $\Om$-invariant complex irreducible characters of $G$ and
$\Irr_\Om(B_0(G))$ to denote the $\Om$-invariant complex irreducible characters
in the principal $p$-block $B_0(G)$. This relates to $\Irr_\cJ(B_0(G))$ thanks to the following lemma.

\begin{lem}
 Let $G$ be a finite group and let $\chi\in\Irr(G)$. Then $\chi$ is
 $\cJ$-invariant if and only if $\chi$ is $\Om(|G|)$-invariant.
\end{lem}

\begin{proof}
Clearly all elements of $\cJ$ restrict to elements of $\Om:=\Om(|G|)$.
Conversely, we show that elements of $\Om$ lift to elements of the same order
in~$\cJ$. It suffices to do this for $\sigma\in\Gal(\QQ_{r^n}/\QQ)$ of
order~$p$, where $r$ runs over primes different from $p$. The hypotheses imply
that $r>2$.
Therefore, for all $m\ge n$, $\Gal(\QQ_{r^m}/\QQ)$ and $\Gal(\QQ_{r^n}/\QQ)$
are cyclic groups with the same $p$-part, with the second group being a factor
group of the first one. In particular, any $p$-element of $\Gal(\QQ_{r^n}/\QQ)$
lifts to a \emph{unique} element of the same order of $\Gal(\QQ_{r^m}/\QQ)$.
Hence $\sigma$ lifts to an automorphism of order $p$ of
$\cup_{m\ge1}\QQ_{r^m}$, which moreover we may choose to act trivially on the
linearly disjoint extensions $\QQ_k/\QQ$ with $k$ prime to~$r$, so to an
element of $\cJ$.
\end{proof}

In particular, we have the following.

\begin{cor}
 Let $H$ be a group of order $n$ and let $\theta\in\Irr(H)$. Suppose that $m$
 is a multiple of $n$. Then $\theta$ is $\Om(n)$-invariant if and only if it
 is $\Om(m)$-invariant.
 \end{cor}

%%%%%%%%%%%%%%%%%%%%%%%%%%%%%%%%%%%%%%%%%%%%%%%%%%%%%%
\subsection{Extensions of It\^o--Michler for almost simple groups}
We consider the following setting: $\bG$ is a simple linear algebraic group
of simply connected type and $F:\bG\to\bG$ is a Steinberg endomorphism with
(finite) group of fixed points $G=\bG^F$. We let $(\bG^*,F)$ be in duality with
$(\bG,F)$ and $G^*:=\bG^{*F}$ (see e.g.~\cite[Def.~1.5.17]{GM20}). For $q$ a 
prime power we denote by $e_p(q)$ the order of $q$ in $\FF_p^\times$.

Recall that a character $\chi$ of a finite group is \emph{$p$-rational} if
its character field $\QQ(\chi)$ is contained in $\QQ_n$ for some $n$ prime
to~$p$.

\begin{lem}   \label{lem:element}
 Let $G=\bG^F$, $G^*=\bG^{*F}$ be as above with $p$ a prime different from the
 underlying characteristic of $\bG$. Assume there is a prime
 $r\not\equiv0,1\pmod p$ with:
 \begin{enumerate}[\rm(1)]
  \item $r$ does not divide $|\bZ(G)|$; and
  \item there is an $r$-element $s\in G^*$ which does not centralise a Sylow
   $p$-subgroup of $G^*$.
 \end{enumerate}
 Then $S=G/\bZ(G)$ has an irreducible $p$-rational and $\Om(|S|)$-invariant
 character of degree divisible by $p$.
\end{lem}

\begin{proof}
Since $r$ does not divide $|\bZ(G)|=|G^*:[G^*,G^*]|$ (see
\cite[Prop.~24.21]{MT}) we have $s\in[G^*,G^*]$ and moreover $\bC_{\bG^*}(s)$
is connected \cite[Prop.~14.20]{MT}. By Lusztig's classification there exists a
unique semisimple character $\chi\in\Irr(G)$ in the Lusztig series of $s$ 
(see \cite[Def.~2.6.9, 2.6.10]{GM20}) and this has $\bZ(G)$ in its kernel, hence
descends to a character of $S$. By Jordan decomposition \cite[Thm~2.6.11]{GM20}
the degree of $\chi$ is divisible by the $p$-part of $|G^*:\bC_{G^*}(s)|$,
hence by $p$. From the explicit definition of $\chi$ as a linear combination of
Deligne--Lusztig characters, the character field $\QQ(\chi)$ is contained in
$\QQ_r$; since $r\ne p$ this shows that $\chi$ is $p$-rational and,
as $r\not\equiv1\pmod p$, also that $\chi$ is $\Om(|G|)$-invariant.
\end{proof}

\begin{prop}   \label{prop:p-rat}
 Let $S$ be a non-abelian simple group and $p>2$ a prime divisor of~$|S|$. If
 all $p$-rational irreducible characters of $S$ which are $\Om(|S|)$-invariant
 have $p'$-degree, then one of the following holds: 
 \begin{enumerate}[\rm(1)]
  \item $S=\tw2\type{B}_2(q^2)$ with $q^2=2^{2f+1}$ and
   $p|(q^2+\eps\sqrt{2}q+1)$, $\eps\in\{\pm1\}$, and all prime divisors of
   $(q^2-1)(q^2-\eps\sqrt{2}q+1)$ are congruent to $1\pmod p$;
  \item $S=\PSL_2(q)$ with $p|(q-\eps)$, and all prime divisors of
   $(q+\eps)/\gcd(2,q+\eps)$ are congruent to~$1\pmod p$;
  \item $S=\PSU_3(q)$ with $q=2^{2f+1}$ and $p=3$, and all prime divisors of
   $(q-1)(q^2-q+1)/3$ are congruent to~$1\pmod3$;
  \item $S=\PSL_4(q)$ with $q=2^{2f+1}$ and $p=3$, and all prime divisors of
   $(q^3-1)(q^2+1)$ are congruent to~$1\pmod3$;
  \item $S=\PSU_4(q)$ with $p|(q-1)$, and all prime divisors $r>2$ of
   $(q^3+1)(q^2+1)$ are congruent to~$1\pmod p$.
 \end{enumerate}
\end{prop}

\begin{proof}
For the sporadic simple groups, the claim is easily checked from the Atlas
\cite{Atl} or using \cite{gap}. Now assume $S=\fA_n$ with $n\ge5$. If $p\ge5$
there exists a $p$-core of size~$n$, that is, an irreducible (and rational)
character of $\fS_n$ of $p$-defect zero. Its restriction to $\fA_n$ has at most
two constituents which hence must be $\Om$-invariant and still of $p$-defect
zero, so $p$-rational.
If $p=3$ then we take for $\chi$ the character of $\fS_n$ labelled by the
partition $(n-2,2), (n-2,1^2)$ of degree~$n(n-3)/2$, $(n-1)(n-2)/2$
respectively, at least one of which is divisible by~$3$.
All of these have irreducible and hence rational restriction to $\fA_n$.
\par
Thus, finally $S$ is of Lie type and $p>2$. If $(S,p)$ is an exception to
the claim, then in particular, all rational irreducible characters of $S$ must
have $p'$-degree. Thus $(S,p)$ appears in the conclusion of \cite[Thm~C]{Gr}.
We discuss these cases, numbered (1)--(7), in turn, striving to exhibit a
suitable prime $r$ such that Lemma~\ref{lem:element} applies. For this,
except for case~(1), we view $S$ as $G/\bZ(G)$ for $G=\bG^F$ for a simple simply
connected group $\bG$ with a Frobenius map $F:\bG\to\bG$. Note that in all
cases $p$ is not the defining characteristic of~$S$, as otherwise the Steinberg
character is as desired (see the proof of \cite[Thm~C]{Gr}).
\begin{enumerate}[(1)]
\item Here $S=\tw2\type{B}_2(q^2)$ with $q^2=2^{2f+1}$ and $p|(q^2+\eps\sqrt{2}q+1)$.
If $r\ne p$ is a prime divisor of $(q^2-1)(q^2-\eps\sqrt{2}q+1)$ with
$r\not\equiv1\pmod p$ then $r$ is as required in Lemma~\ref{lem:element}.
\item Here $S=\PSL_n(q)$ and $n<p|(q-1)$. We may assume $n>2$ since the
case $n=2$ is treated in (4) below. Then we have
$$(q^{n-1}-1)/(q-1)\equiv n-1\pmod{q-1}\equiv n-1\pmod p\not\equiv 0,1\pmod p,$$
hence $(q^{n-1}-1)/(q-1)$ has a prime divisor $r\not\equiv0,1\pmod p$. Since
$r$ divides $(q^{n-1}-1)/(q-1)$ we have $r$ is prime to $|\bZ(G)|=\gcd(n,q-1)$,
with $G=\SL_n(q)$.
Let $s$ be a (semisimple) element of $G^*=\PGL_n(q)$ generating the Sylow
$r$-subgroup of a (cyclic) maximal torus of order $q^{n-1}-1$. As $s$ is not
contained in a maximally split torus of $\PGL_n(q)$ (whose exponent is $q-1$,
while by the choice of $r$, the order of $s$ is larger than $(q-1)_r$), $s$
does not centralise a maximally split torus of $\PGL_n(q)$, hence does not
centralise a Sylow $p$-subgroup of $G^*$. Thus all assumptions of
Lemma~\ref{lem:element} are satisfied and this case is not an exception.
\item Here $S=\PSU_n(q)$ and $2<n<p|(q+1)$. Now
$$(q^{n-1}-(-1)^{n-1})/(q+1)\equiv n-1\pmod{q+1}\equiv n-1\pmod p
  \not\equiv 0,1\pmod p,$$
hence $(q^{n-1}-(-1)^{n-1})/(q+1)$ has a prime divisor $r\not\equiv0,1\pmod p$,
and as before, $r$ is prime to $|\bZ(G)|=\gcd(n,q+1)$, with $G=\SU_n(q)$. The
rest of the argument is now entirely analogous to the previous case to show
that Lemma~\ref{lem:element} applies.
\item Here $S=\PSL_2(q)$, $S=\PSL_4(q)$ with $e_p(q)=2$, or $S=\PSL_3(q)$
with $p=3|(q-1)$. In the latter case, $q+1\equiv2\pmod3$, so $q+1$ has a
prime divisor $r\equiv2\pmod3$. Let $s$ be an element of order $r$ in the dual
group $G^*=\PGL_3(q)$. Since $r$ divides $q+1$, the centraliser of $s$ in $G^*$
does not contain a Sylow 3-subgroup of $G^*$, and $r$ does not divide
$|\bZ(\SL_3(q))|$. Thus, Lemma~\ref{lem:element} applies.
Now assume $S=\PSL_2(q)$. If $p|(q-\eps)$ and $q+\eps$ has a prime divisor
$r>2$ with $r\not\equiv1\pmod p$ then we are done by Lemma~\ref{lem:element}.
Similarly, if $2|(q+\eps)/2$ then the characters of $S$ of degree $(q-\eps)/2$
have values in $\QQ(\sqrt{\pm q})$ (see \cite[Tab.~2.6]{GM20}), so
are $p$-rational and $\Om$-invariant.
Finally, let $S=\PSL_4(q)$ with $2<p|(q+1)$. 
If $p>3$ then let $r$ be any prime divisor of $q-1\equiv-2\not\equiv1\pmod p$
which is not congruent to~$1\pmod p$ and $s$ an element of maximal $r$-power
order in $G^*$. If $p=3$ and $q$ is odd, let $s\ne 1$ be a 2-element in
$\PSL_4(q)=[G^*,G^*]$ not centralising a Sylow 3-subgroup of $G^*$, and we
may again conclude. Finally, when $p=3$ and $q$ is even (and hence an odd
2-power), then whenever $r$ is any prime divisor of $(q^3-1)(q^2+1)$ not
congruent to~$1\pmod 3$, we can again apply Lemma~\ref{lem:element}.
\item Here $S=\PSU_4(q)$ with $p|(q-1)$ or $S=\PSU_3(q)$ with $p=3|(q+1)$.
In the latter case, if $q$ is odd then take for $s$ a 2-element of $G^*$ of
maximal possible order (and hence not centralising a Sylow $3$-subgroup), and
apply Lemma~\ref{lem:element}. If
$q$ is even, and there is a prime $r$ not congruent to~$1\pmod3$ dividing
$(q-1)(q^2-q+1)/3$, this will do.  Finally, if $S=\PSU_4(q)$ with $p|(q-1)$,
we are done by Lemma~\ref{lem:element} if there is a prime divisor $r>2$ of
$(q^3+1)(q^2+1)$ not congruent to $1\pmod p$.
\item Here $S=\PSO_{2n}^+(q)$ with $n\in\{5,7\}$ and $p|(q+1)$. If $p=3$ any
(rational) unipotent character of degree divisible by~$q^2-q+1$, hence by~$p$,
is as desired. So assume $p>3$. Then there exists a prime
divisor $r\ne p$ of $q^2-q+1\equiv3\pmod p$ not congruent to~1 modulo~$p$.
Let $s$ be an element generating a Sylow $r$-subgroup of a torus of order
$q^3+1$ in $G^*$, where $G=\Spin_{2n}^+(q)$. Since $r\ne2$ this satisfies the
assumptions of Lemma~\ref{lem:element}.
\item Here $S=\PSO_{2n}^-(q)$ with
$$(n,e_p(q))\in\{(4,1),(4,2),(4,4),(5,1),(6,1),(6,2),(7,1),(8,1),(8,2)\}.$$
First assume $e_p(q)=1$, so $p|(q-1)$. If $p=3$ we use a (rational) unipotent
character of degree divisible by $q^2+q+1$, so by~3. If $p\ne3$ there exists a
prime divisor $r\ne p$ of $q^2+q+1\equiv3\pmod p$ not congruent to~1 modulo~$p$,
and any $r$-element generating a torus of order $q^3-1$ in the dual
group $G^*$ of $G=\Spin_{2n}^-(q)$ satisfies the assumptions of
Lemma~\ref{lem:element}. Next assume $e_p(q)=2$. Again, when $p=3$ there is a
suitable unipotent character of degree divisible by~$q^2-q+1$. Otherwise, let
$r$ be a prime divisor of $q^2-q+1\equiv3\pmod p$ not congruent to~1
modulo~$p$, and argue as before. So finally
$S=\PSO_{8}^-(q)$ and $2\ne p|(q^2+1)$. With $r$ a prime divisor of
$q^2+q+1\equiv q\not\equiv1\pmod p$ not congruent to~1 modulo~$p$, we are done
as before.\qedhere
\end{enumerate}
\end{proof}

\begin{exmp}   \label{ex:PSL}
 The second case listed in Proposition~\ref{prop:p-rat} is a true exception,
 for example the groups $S=\PSL_2(q)$ with
 $$(p,q)\in
   \{(3,8),(5,61),(7,421),(11,397),(13,157),(17,613),(19,457),(23,277)\},$$
 and similarly the first case is a true exception for example when
 $S=\tw2\type{B}_2(8)$ with $p=5$. We have not tried to determine whether there
 are cases with arbitrarily large $p$.
\end{exmp}

\begin{exmp}   \label{example:J1}
 The sporadic simple group $J_1$ has only three irreducible characters of
 degree divisible by~3; they are cyclically permuted by the Galois automorphism
 of order~9 of $\QQ_{19}/\QQ$ (see \cite{Atl}), but fixed by the elements of
 order~3 in $\Om(19)$.
\end{exmp}

\begin{thm}   \label{thm:as}
 Let $A$ be an almost simple group with socle $S$ such that $A/S$ is a
 $p$-group for a prime $p$ dividing $|A|$. If $A$ has no $\Om(|A|)$-invariant
 irreducible character of degree divisible by $p$, then $p\ne2$ and
 $(S,p)$ is as in (1)--(5) of Proposition~\ref{prop:p-rat}.
\end{thm}

\begin{proof}
First assume $p$ is odd.
Assume there is a $p$-rational and $\Om$-invariant character $\chi\in\Irr(S)$ of
degree divisible by~$p$. Then by \cite[Thm~6.30]{isa}, $\chi$ has a unique
$p$-rational extension $\tilde\chi$ to its inertia group $I_A(\chi)$, which
then, by uniqueness, must also be $\Om$-invariant. Then $\tilde\chi^A$ is a
character of $A$ of degree divisible by~$p$ as desired. Hence, we arrive at
$(S,p)$ being one of the exceptions in Proposition~\ref{prop:p-rat}.
\par
This leaves the case when $p$ does not divide $|S|$. Then necessarily $A/S$
is generated by field automorphisms (see e.g.~\cite{Atl}), say
$A/S=\langle\gamma\rangle$ with $|\gamma|=p^a$. Let $r$ be a prime dividing
$|S|$ but not $|\bC_S(\gamma^{p^{a-1}})|$, which exists by Zsigmondy's theorem,
and let $s\in[G^*,G^*]$ be an $r$-element. Then the semisimple character
$\chi_s$ of $G$, considered as character of $S$, has values in $\QQ_{r^b}$ for
some $b\ge0$. Since $r$ does not divide $|\bC_G(\gamma^{p^{a-1}})|$, $\chi_s$
lies in an orbit of length $|\gamma|$ under $A$, so the induced character
$\chi^A$ is irreducible. By the choice of $s$, the values of $\chi^A$ lie in a
subfield of index $|\gamma|=p^a>1$ of $\QQ_{r^b}$, which is $\Om$-invariant.

Now assume $p=2$. By \cite[Thm~C]{TT22}, there is a rational-valued character
in $\Irr(A)$ of even degree unless $S=\PSL_2(q)$ with $q=3^f$ and $f\ge3$ is
odd. Consider the latter case. Since $A/S$ is a $2$-group, we have
$A\in\{S, \wt{S}\}$, where $\wt{S}:=\PGL_2(q)$. Let $r$ be a prime divisor of
$(q-1)/2$ and $\chi$ an irreducible Deligne--Lusztig character of $\SL_2(q)$
labelled by an element in $\PGL_2(q)$ of order $r$, so of even degree~$q+1$.
Then $\chi$ has values in $\QQ_r$ by \cite[Tab.~2.6]{GM20} and it has
$\bZ(\SL_2(q))$ in its kernel as $r$ is odd. The unique involution in
$\Gal(\QQ_r/\QQ)$ acts as complex conjugation, but $\chi$ is real by loc.~cit.
Since $\QQ(\chi)\le\QQ_r$, then $\chi$ is $\Om(|S|)$-invariant. Moreover,
$\chi$ is in fact the restriction of an irreducible character in $\Irr(\wt{S})$
with the same field of values, so the proof is complete.
\end{proof}

The following result on quasi-simple groups will be needed in the proof of
Theorem~\ref{main}:

\begin{thm}   \label{thm:quasisimples}
 Let $p$ be a prime number. Let $S$ be a non-abelian simple group with abelian
 Sylow $p$-subgroups. If $G$ is quasi-simple with $G/\bZ(G)=S$ and
 $|\bZ(G)|=p$, then there exists a faithful
 $\chi\in\Irr_{\Om(|G|)}(B_0(G))$ (and hence of degree divisible by $p$).
\end{thm}

\begin{proof}
First assume $p=2$. The non-abelian simple groups with abelian Sylow
2-subgroups are $\PSL_2(q)$ with $q\equiv3,5\pmod8$ or $q$ even, the Ree groups,
and $J_1$. Of these, only $\PSL_2(q)$ with $q\equiv3,5\pmod8$ have covering
groups $G=\SL_2(q)$ as in the claim. Now $\SL_2(q)$ has a rational faithful
irreducible Deligne--Lusztig character of degree $q-1$ (if $q\equiv3\pmod8$)
respectively $q+1$ (if $q\equiv5\pmod8$) labelled by an element of order~4 in
the dual group $\PGL_2(q)$, hence lying in the principal 2-block (e.g.~by
\cite[Prop.~6]{En00}).
\par
So now assume $p\ge3$. The non-abelian simple groups with order of the Schur
multiplier divisible by an odd prime $p$ are $\PSL_n(\eps q)$ with
$p|\gcd(n,q-\eps)$, and for $p=3$ groups of type $E_6(\eps q)$, $\fA_6$,
$\fA_7$, and some of the sporadic simple groups. Among those, the only ones
with abelian Sylow $p$-subgroups occur when $p=3$, namely $\PSL_3(q)$ with
$q\equiv4,7\pmod9$, $\PSU_3(q)$ with $q\equiv2,5\pmod9$, $\fA_6$, $\fA_7$,
$M_{22}$ and $ON$. Here, the groups $\PSL_3(\eps q)$ have faithful
characters of degree $(q^2-1)(q-\eps)$ with values in $\QQ_3$, labelled by an
element of order~3 in $\PGL_3(\eps q)$, so lying in the principal 3-block.
Using \cite{gap}, the groups $3.\fA_6$, $3.\fA_7$, $3.M_{22}$, and $3.ON$ are
seen to have faithful characters of degrees 6, 15, 231, 495 respectively lying
in the principal 3-block, with values in $\QQ_3$. This completes the proof.
\end{proof}

%%%%%%%%%%%%%%%%%%%%%%%%%%%%%%%%%%%%%%%%%%%%%%%%%%%%%%
\subsection{Some observations on extending unipotent characters}
We continue to consider connected reductive linear algebraic groups $\bG$ with
a Steinberg endomorphism $F:\bG\to\bG$.
The next result will be useful when dealing with graph automorphisms.

\begin{lem}[Digne--Michel]   \label{lem:unipextgraph}
 Let $\bG$ be of type $\type{A}_{n-1}$, $\type{E}_6$, or $\type{D}_n$ with
 $F$ inducing an $\FF_q$-structure, and $G=\bG^F$. Let $\tau$ be a non-trivial
 graph automorphism of $G$. Then any $\tau$-invariant unipotent character
 $\chi$ lying in the principal series of $G$ extends to a rational-valued
 character $\tchi$ of $G\langle\tau\rangle$, unless either $G=\type{E}_6(q)$
 and $\chi\in\{\phi_{64,4},\phi_{64,13}\}$, or $G=\type{A}_{n-1}(q)$ and $\chi$
 is labelled by a partition $\la=(\la_1,\ldots,\la_r)$ of $n$ with
 \[\sum_i {{\la_i}\choose{2}}-\sum_j {{\la_j'}\choose{2}} + {{n}\choose{2}}
   \equiv1\pmod 2,\]
 where $\la'=(\la_1',\ldots,\la_s')$ is the partition conjugate to $\la$.
 In the latter cases, $\QQ(\tchi)=\QQ(\sqrt{q})$.
\end{lem}

\begin{proof}
This follows from \cite[Thm~II.3.3]{DM85} (see also \cite[Thm~3]{DM24}).  Note
that all unipotent characters under consideration are rational-valued, by
\cite[Thm~2.9]{BC72}. 
\end{proof}

Given a normal subgroup $N$ of a group $G$ and $\theta\in\Irr(N)$, we will
write $I_G(\theta)$ or $G_\theta$ to denote the inertia group of $\theta$
in~$G$. In the case of field automorphisms, the following is also useful:

\begin{lem}   \label{lem:unipprinc}
 Let $\bG$ be simply connected with a Frobenius map $F$ such that $S=G/\bZ(G)$
 is simple, where $G:=\bG^F$. Let
 $\chi\in\Irr(S)$ be the deflation of a principal series unipotent character.
 If $S\leq A\leq\Aut(S)$ is such that $A$ is generated by inner-diagonal and
 field automorphisms, then $\chi$ has an extension $\hat\chi$ to $I_A(\chi)$
 with $\QQ(\hat\chi)\le\QQ(\chi)$.
\end{lem}

\begin{proof}
By \cite[Thm~2.4]{Ma08}, every unipotent character extends to its stabilizer in
$\Aut(S)$. The claim then follows directly from \cite[Prop.~2.6]{RSF22} (and
its generalization \cite[Prop.~8.7]{Jo22}), as in the notation of loc.~cit.,
the cuspidal character $\delta=1$ extends trivially to $\hat\delta=1$.
%If $S=\tw{2}\type{B}_2(q^2)$ or $\tw{2}\type{G}_2(q^2)$, then $\chi$ is the
%Steinberg character or trivial character, and the result follows by\cite{schmid}.
\end{proof}

Using the above lemma, the following is proved in Lemma 4.4 of \cite{MMSV}:

\begin{lem}   \label{lem:unipext}
 Let $p$ be a prime and $\bG$ be simply connected defined in characteristic
 distinct from~$p$, with a Steinberg map $F$ such that $S=G/\bZ(G)$ is simple,
 where $G:=\bG^F$. 
 Then there exists a non-trivial, rational-valued unipotent character
 $\chi\in\Irr_{p'}(B_0(S))$ such that if $S\leq A\leq\Aut(S)$ and $A$ is
 generated by inner-diagonal and field automorphisms, then there is an
 extension $\hat\chi$ of $\chi$ to $I_A(\chi)$ that is rational-valued.
\end{lem}

We will also use the following in the case $p=2$.

\begin{prop}   \label{prop:evenrat}
 Let $S$ be a simple group of Lie type defined in odd characteristic, different
 from
 types $\type{A}_1$, $\tw2\type{A}_2$ and $\tw2\type{G}_2$. Then the principal
 $2$-block $B_0(S)$ contains a rational unipotent principal series character
 $\chi$ of even degree such that if $S\leq A\leq\Aut(S)$ and $A$ is generated
 by inner-diagonal and field automorphisms, then $\chi$ has a rational
 extension $\hat\chi$ to $I_A(\chi)$.
\end{prop}

\begin{proof}
Suppose first that $S$ is of classical type. Then all unipotent characters of
$S$ are rational (see \cite[Cor.~4.4.24]{GM20}) and lie in the principal
2-block (see \cite[Prop.~6]{En00}). Since the Weyl group of $S$ has an
irreducible character of
even degree unless we are in the excluded cases, the claim follows by the
arguments in \cite[2.3]{KM17}. Now assume $S$ is of exceptional type. Then by
\cite[Thm~A]{En00} all principal series unipotent characters lie in $B_0(S)$.
It is easy to check with \cite[Cor.~4.5.6]{GM20} that there is a rational one
of even degree.   \par
The second claim now follows from Lemma \ref{lem:unipprinc}. 
\end{proof}

%%%%%%%%%%%%%%%%%%%%%%%%%%%%%%%%%%%%%%%%%%%%%%%%%%%%%%
\subsection{Existence of certain character extensions}
We work towards a proof of Conjecture \ref{galcon} and begin with some
additional results about extensions.

When $S$ is a simple group of Lie type, we write $\Aut(S)=\wt{S}\rtimes D$,
where $\wt{S}$ is the group of inner-diagonal automorphisms and $D$ is an
appropriate group of graph-field automorphisms, see \cite[Thms~2.5.12
and~2.5.14]{GLS}.

\begin{thm}   \label{thm:ext}
 Let $A$ be an almost simple group with socle $S$. Suppose that $A/S$ is a
 $p$-group.  Then there exists $1_S\neq\vhi\in\Irr_\Om(S)$ that extends to an
 $\Om$-invariant character $\hat\vhi$ of $I_A(\vhi)$. Furthermore, if $p$
 divides $|S|$ we can take $\hat\vhi$ to lie in the principal $p$-block
 of $I_A(\vhi)$. In this case, $\hat\vhi\in\Irr(B_0(I_A(\vhi)))$ can be chosen
 to be rational-valued, except possibly when $S=\PSL_2(p^f)$, or when $p=2$ and
 $S=\PSL_3(\pm 2^f)\}$ or $A=S=\tw2\type{B}_2(q^2)$.
\end{thm}

\begin{proof}
We see from \cite[Lemma~4.1]{HSTV} that there is a non-trivial rational-valued
character $\vhi\in\Irr(S)$ that extends to a rational-valued character of
$\Aut(S)$, so the first statement holds. 

For the second statement, since $A/S$ is a $p$-group and hence $B_0(I_A(\vhi))$
is the unique $p$-block of $I_A(\vhi)$ above $B_0(S)$, it suffices to show that
a $\vhi$ satisfying the first statement can be chosen in $B_0(S)$. For the
third statement, we wish to ensure that further $\hat\vhi$ is rational-valued,
when $S$ is not one of the listed exceptions.

For the sporadic groups and the alternating group $\alt_6$, this may be checked
in GAP \cite{gap}. Next, suppose that $S=\alt_n$ with $n\geq 5$, $n\ne6$. Then
$A\in\{\alt_n, \sym_n\}$.
The character of $\sym_n$ corresponding to any non-self-conjugate partition
of~$n$ with $p$-core $(r)$, where $n=ap+r$ with $0\leq r<p$, will lie in
$B_0(\sym_n)$ and restrict irreducibly to $B_0(\alt_n)$, giving a rational
character of $S$ that extends to a rational character in $B_0(A)$. More
specifically, we may take the partition $(n-1,1)$ if $r=0$; $(ap-1,2)$ if
$r=1$; and $(r,1^{ap})$ if $r>1$.

Now let $S$ be a simple group of Lie type. First, assume that $p$ is not the
defining characteristic for $S$. Suppose that $p\in\{2,3\}$. Then the Steinberg
character $\St_S$ lies in $B_0(S)$, using \cite[Thm~A]{En00} and
\cite[Thm~6.6]{Ma07}, and extends to a
rational-valued character of $\Aut(S)$ by \cite{schmid}.
If instead $p\geq 5$, note that $A$ is generated by inner-diagonal and field
automorphisms, since $A/S$ is a $p$-group. Then the result follows from
Lemma~\ref{lem:unipext}.

Now assume that $S$ is a group of Lie type defined in characteristic $p$.
We have $\Irr(B_0(S))=\Irr(S)\setminus\{\St_S\}$ by \cite{Hum}, so
it suffices to know that there is a character $\vhi\not\in\{1_S, \St_S\}$ that
extends to an $\Om$-invariant (in fact rational, except in the stated
exceptions) character of $I_A(\vhi)$. Recall that $\Aut(S)=\wt{S}\rtimes D$. We
have $p\nmid |\wt{S}/S|$ and all Sylow subgroups of $D$ are abelian, and hence
the Sylow $p$-subgroups of $\Aut(S)/S$ are abelian. In particular, $A/S$ is
abelian.

If $S=\tw{2}\type{F}_4(q^2)$ or $\tw{2}\type{B}_2(q^2)$, then $p=2$ and
$|\Aut(S)/S|$ is odd, so in this case $A=S$. When $S=\tw{2}\type{B}_2(q^2)$ we
take for $\vhi$ one of the two cuspidal unipotent characters, with values in
$\QQ_4$, and when $S=\tw{2}\type{F}_4(q^2)$ any non-trivial rational unipotent
character apart from the Steinberg character. Observe that for
$\tw{2}\type{B}_2(q^2)$ there is no suitable rational character. 

If $S=\tw{2}\type{G}_2(q^2)$, then $p=3$ and we may take $\vhi$ to be the
semisimple character of degree $q^4-q^2+1$ labelled by the class of involutions
in $G^*\cong S$. Then $\vhi$ is rational-valued, and note that $3\nmid\vhi(1)$.
Since $A/S$ is a $3$-group, $\vhi$ extends to a rational-valued character
of~$A$ by \cite[Cors~6.2 and 6.4]{N18}.
We therefore assume that $S$ is not a Suzuki or Ree group.

First suppose that $S=\PSL_n(\eps q)$ with $n\geq 4$, $q=p^f$, and
$\eps\in\{\pm1\}$. Write $D=\langle \tau, F_p\rangle$, where $\tau$ is a graph
automorphism of order $2$ and $F_p$ is a generating field automorphism of
order~$f$ if $\eps=1$, respectively~$2f$ if $\eps=-1$. 
Let $\la$ denote the partition $(\frac{n-4}{2}+2,2,1^{\frac{n-4}{2}})$ if
$n\equiv0\pmod4$, the partition $(\frac{n-1}{2},1^{\frac{n-1}{2}})$ if
$n\equiv1\pmod4$, and the partition $(n-2,2)$ if $n\equiv2,3\pmod4$.
Then $\la$ does not satisfies the condition in Lemma~\ref{lem:unipextgraph}.
Hence the
corresponding unipotent character $\chi_\la\in\Irr(\wt{S})$ extends to a
rational-valued character of $\wt{S}\langle\tau\rangle$. As $\chi_\la$ lies
in the principal series, it further extends to a rational-valued character of
$\wt{S}\langle F_p\rangle$, by Lemma~\ref{lem:unipprinc}. Note that if
$p\geq 3$ or $f$ is odd, then $A\leq \wt{S}\langle F_p\rangle$ or
$A\leq \wt{S}\langle\tau\rangle$, since $A/S$ is a $p$-group. Then
$\vhi:=(\chi_\la)_S$ extends to a rational-valued character of $I_A(\vhi)$ in
these cases. 
 
So now assume $p=2$ and $f$ is even. Since $\Aut(S)/\wt{S}$ is abelian and any
unipotent character extends to $\Aut(S)$ by \cite[Thms~2.4, 2.5]{Ma08}, we have
any character of $\Aut(S)$ lying above $\chi_\la$ is an extension by
Gallagher's theorem.  From above, any extension of $\chi_\la$ to
$\wt S\langle\tau\rangle$ is rational. Choosing $\hat\vhi$ an extension of
$\chi_\la$ to $\Aut(S)$ to be such that $\hat\vhi|_{\wt{S}\langle F_p\rangle}$
is the rational extension of $\vhi$ guaranteed by Lemma~\ref{lem:unipprinc},
we obtain $\hat\vhi|_A$ is rational-valued.
 
Next, consider the case $S=\PSL_3(\eps q)$ with $q=p^f$. We remark first that
if $\eps=1$ and either $f$ is even or $p$ is odd, the same arguments as above
but with $\la=(2,1)$ yield a unipotent character $\vhi\not\in\{1_S, \St_S\}$
that extends to a rational character of $I_A(\vhi)$. 

Let $\wt{G}:=\GL_3(\eps q)$ and let $G:=\SL_3(\eps q)$. 
Assume $p$ is odd. Let $s=\diag(-1,-1,1)\in G\cong[\wt G^*,\wt G^*]$. Then
the semisimple character $\chi_s$ of $\wt{G}$ is trivial on $\bZ(\wt{G})$,
restricts irreducibly to $G$, and has degree prime to $p$. Further, $\chi_s$ is
rational-valued. Then taking $\vhi\in\Irr(S)$ to be the deflation of
$\chi_s|_G$ and applying \cite[Cors~6.2 and~6.4]{N18}, we see $\vhi$ extends
to a rational-valued character of $I_A(\vhi)$.

Now assume $p=2$. Let $s=\diag(\ze,\ze^{-1},1)\in G\cong[\wt G^*,\wt G^*]$ with
$3<|\ze|$ a prime power dividing $q-\eps$. Note that this is possible since
$q\neq 2$ in the case $\eps=-1$
as $S$ is assumed simple, and $q\neq 4$ in the case $\eps=1$ since $f$ is odd.
Then the semisimple character $\chi_s\in\Irr(\wt{G})$ restricts irreducibly
to~$G$, is trivial on $\bZ(\wt{G})$, and has odd degree. Further,
$\ze^\sigma\in\{\ze,\ze^{-1}\}$ for any Galois automorphism $\sigma$ of
order~$2$ since $(\ZZ/|\ze|\ZZ)^\times$ is cyclic.
%$\ze^\sigma=\ze^k$ for some $k$ relatively prime to $|\ze|$, and
%hence $\ze=\ze^{\sigma^2}=\ze^{k^2}$, so the odd prime power $|\ze|$ divides
%$k^2-1$, and hence either $k-1$ or $k+1$.
Then $\chi_s^\sigma=\chi_s$ for any
$\sigma\in\Om$. Let $\vhi\in\Irr(S)$ be the deflation of the restriction of
$\chi_s$ to $G$. %Then by \cite[Prop.~6.1]{Ru21}, 
Then again by \cite[Cors~6.2 and~6.4]{N18}, we have 
an extension of $\vhi$ to $I_A(\vhi)$ that is also $\Om$-invariant.
 
Now consider the case $S=\PSL_2(q)$ with $q=p^f$. If $q\leq 9$, then we can
check directly in GAP or use the well-known isomorphisms with alternating
groups to see the statement holds. So, assume $q\geq 11$. If $p$ is odd, then
by \cite[Tab.~2.6]{GM20} there are two characters of $S$ of degree
$\frac{q+\eta}{2}$, where $\eta\in\{\pm1\}$ with $q\equiv\eta\pmod 4$, which
are $\Om$-stable and have $p'$-degree, so again possess an $\Om$-invariant
extension to
$I_A(\vhi)$ by \cite[Cors~6.2 and 6.4]{N18}. We may therefore assume that
$p=2$. Let $\wt{G}:=\GL_2(q)$ and note that $S=\SL_2(q)$. Let $s\in\wt{G}^*$ be
an element with eigenvalues $\{\ze,\ze^{-1}\}$, where $|\ze|$ is a
prime dividing $q+1$. If $\sigma\in\Om$ has order~$2$, then as above,
$\ze^\sigma\in\{\ze,\ze^{-1}\}$, so the corresponding semisimple
character $\chi_s\in\Irr(\wt G)$ is $\sigma$-invariant, whence
$\chi_s\in\Irr_\Om(\wt{G})$. Further, $\chi_s$ has odd degree and restricts
irreducibly to $S$. Then letting $\vhi$ denote this restricted character, we
have $\vhi\in\Irr_{2',\Om}(B_0(S))$ since $\vhi\neq\St_S$. Again applying
\cite[Cors~6.2 and 6.4]{N18}, $\vhi$ has an $\Om$-invariant extension
to~$I_A(\vhi)$. 

We may now assume that $S=G/\bZ(G)$ where $G$ is not of type $\type{A}$,
$\tw{2}\type{A}$, or Suzuki or Ree type.  The principal series unipotent
characters are rational-valued except for a small number of exceptions for
$S=\type{E}_7(q)$ and $\type{E}_8(q)$, by \cite[Thm~2.9]{BC72}. Further, in our
remaining cases, there is always a principal series unipotent character
$\vhi\not\in \{1_S, \St_S\}$ (and hence in $B_0(S)$) that is rational-valued,
distinct from the exceptions for $\type{E}_6(q)$ in
Lemma~\ref{lem:unipextgraph}, and of degree divisible by~$p$. (The principal
series unipotent characters are described in \cite[Sec.~13.8, 13.9]{Ca85}.) In
the case of $\type{B}_2(2^{n})$ with $n\geq 2$ and $\type{F}_4(2^n)$ with
$n\geq 1$, we may take $\varphi$ to be the character indexed by the symbol
${1, 2}\choose {0}$, respectively the character $\phi_{8,3}'$ in the notation
of \cite[Sec.~13.9]{Ca85}, which is stable under field automorphisms but moved
by the exceptional graph automorphism by \cite[Thm.~2.5]{Ma08}. That is, in the
latter cases, $I_A(\varphi)$ is generated by $S$ and field automorphisms.
Then using Lemmas~\ref{lem:unipprinc} and~\ref{lem:unipextgraph}, we can argue
analogously to before to see the statement holds.
\end{proof}

The following will be useful toward proving Conjecture \ref{galcon} for certain
almost simple groups (Theorem \ref{thm:as2} below).

\begin{lem}   \label{lem:oddext}
 Let $p$ be a prime and let $A$ be an almost simple group such that
 $\bO^{p'}(A)=A$ with socle $S$ a simple group of Lie type. Assume $p\geq 5$ if
 $S=\type{D}_4(q)$. Assume that one of the following holds:
 \begin{enumerate}[\rm(1)]
  \item $p$ is odd and $\Irr_\Om(B_0(\wt{S}))$ contains a $p$-rational
   character of degree divisible by $p$ that restricts irreducibly to $S$;
  \item $p$ is odd and $\Irr_\Om(B_0(A\cap\wt{S}))$ contains a $p$-rational
   character not invariant under $A$;
  \item $\Irr_\Om(B_0(A\cap\wt{S}))$ contains a character of $p'$-degree not
   invariant under $A$.
 \end{enumerate} 
 Then there is a character in $\Irr_\Om(B_0(A))$ of degree divisible by $p$.
\end{lem}

\begin{proof}
Since $\bO^{p'}(A)=A$ and either $D$ is abelian or $S=\type{D}_4(q)$ and $D$ is
$\sym_3\rtimes C$ with $C$ cyclic, we have
$A/(A\cap\wt{S})\cong A\wt{S}/\wt{S}$ is a $p$-group.  Then $B_0(A\wt{S})$ is
the unique $p$-block covering $B_0(\wt{S})$ and $B_0(A)$ is the unique block
covering $B_0(A\cap\wt S)$.
\par\noindent
(1) If $\Irr_\Om(B_0(\wt{S}))$ contains a $p$-rational character of degree
divisible by $p$ that restricts irreducibly to $S$, then let $\vhi$ be its
restriction to $A\cap \wt S$. Then $\vhi\in\Irr_\Om(B_0(A\cap\wt{S}))$ is
$p$-rational and has degree divisible by $p$.  By \cite[Thm~6.30]{isa}, if
$p$ is odd then $\vhi$ extends to a unique $p$-rational character $\hat\vhi$ of
$I_A(\vhi)$, which is therefore in $\Irr_\Om(B_0(I_{A}(\vhi)))$. Then the
induced character $\hat\vhi^A$ has degree divisible by $p$ and lies in
$\Irr_\Om(B_0(A))$.
\par\noindent
(2) Now assume $\Irr_\Om(B_0(A\cap\wt{S}))$ contains a $p$-rational character
$\vhi$ not invariant under $A$. Then $[A:I_{A\cap\wt{S}}(\vhi)]$ is a power
of~$p$. As before, if $p$ is odd then $\vhi$ has a unique $p$-rational
extension $\hat\vhi$ to $I_{A}(\vhi)$, which must lie in $B_0(I_A(\vhi))$ and
be $\Om$-invariant. Then $\hat\vhi^{A}$ has degree divisible by $p$ and lies in
$\Irr_\Om(B_0(A))$.
\par\noindent
(3) If instead $\vhi\in\Irr_{p',\Om}(B_0(A\cap\wt{S}))$ is not invariant
under~$A$, we may argue similarly, using \cite[Cors~6.2 and 6.4]{N18} in place
of~\cite[Thm~6.30]{isa}.
\end{proof}

\begin{rem}
Note that if $S=\type{D}_4(q)$ and $A/(A\cap\wt{S})$ is a $p$-group, then the
conclusion of Lemma~\ref{lem:oddext}(1, 2) still holds when $p=3$ and that of
Lemma~\ref{lem:oddext}(3) when $p\in\{2,3\}$.
\end{rem}

%%%%%%%%%%%%%%%%%%%%%%%%%%%%%%%%%%%%%%%%%%%%%%%%%%%%%%
\subsection{Extension of Brauer's height zero conjecture for almost simple groups}
If $p$ is a prime, we denote by $\bO^{p'}(G)$ the smallest normal subgroup of
the group $G$ with factor group of $p'$-order. The following is
Conjecture~\ref{galcon} for certain almost simple groups. It will be used in
the proof of Theorem~\ref{main}.

\begin{thm}   \label{thm:as2}
 Let $A$ be an almost simple group such that $\bO^{p'}(A)=A$. Then $A$ has
 abelian Sylow $p$-subgroups if and only if all the characters in
 $\Irr_\Om(B_0(A))$ have $p'$-degree.
\end{thm}

\begin{proof}
The ``only if" direction follows from that direction of the ordinary height
zero conjecture,
proved in \cite{KM}. For the converse, assume that $A$ has non-abelian Sylow
$p$-subgroups. We claim that $\Irr_\Om(B_0(A))$ contains a character with
degree divisible by $p$. Let $S$ be the simple socle of $A$. We deal with the
various possibilities for $S$ in the subsequent
Propositions~\ref{prop:as2 not Lie}--\ref{prop:as2 cross}.
\end{proof}

\begin{prop}   \label{prop:as2 not Lie}
 The conclusion of Theorem~\ref{thm:as2} holds for $S$ not of Lie type.
\end{prop}

\begin{proof}
For $S=\fA_n$ with $n\ge5$ the claim for $p=2$ is shown in the proof of
\cite[Prop.~2.5]{MN}. For $p>2$ the principal $p$-block of $\fS_n$ contains a
character $\chi$ of degree divisible by~$p$ (if $\fS_n$ has non-abelian Sylow
$p$-subgroups) by Brauer's height zero conjecture \cite{MNST}, and this is
rational.
Furthermore, the constituents of its restriction to $\fA_n$ are $\Om$-invariant
and still of degree divisible by~$p$. If $S$ is sporadic or the Tits simple
group, it can be checked with \cite{gap} that the only cases when $|A|$ is
divisible by $p^3$ but $A$ does not possess a rational character in the
principal $p$-block of degree divisible by~$p$ is when $A=S=J_1$ with $p=2$,
or $S= ON$ with $p=3$. But in either case, Sylow $p$-subgroups of $A$ are
abelian.
\end{proof}

For the groups of Lie type, as before, we write $\Aut(S)=\wt{S}\rtimes D$.
Since $\bO^{p'}(A)=A$, it follows that any field automorphism in $A$ has
$p$-power order, and further any graph-field automorphism in $A$ has $p$-power
order except possibly if $S=\type{D}_4(q)$.

\begin{prop}   \label{prop:as2 defchar}
 The conclusion of Theorem~\ref{thm:as2} holds for $S$ of Lie type in
 characteristic~$p$.
\end{prop}

\begin{proof}
Assume $S$ is simple of Lie type defined in characteristic $p$. Recall again
that $\Irr(B_0(S))=\Irr(S)\setminus\{\St_S\}$. If $S=\tw{2}\type{F}_4(q^2)$ or
$\tw{2}\type{B}_2(q^2)$ with $q^2=2^{2m+1}$ and $p=2$, then the condition
$\bO^{p'}(A)=A$ means that $A=S$. In the case of $S=\tw{2}\type{F}_4(q^2)$,
there is a unique unipotent character of degree $q^2(q^4-q^2+1)(q^8-q^4+1)$,
which is therefore rational-valued and has even degree. When
$S=\tw{2}\type{B}_2(q^2)$, there are two cuspidal unipotent characters of even
degree, which take values in $\QQ(\sqrt{-1})$, and therefore are $\Om$-stable
since $p=2$.

Next suppose $S=\tw{2}\type{G}_2(q^2)$ with $q=3^{2m+1}$ and $p=3$. Then
$\Out(S)$ is cyclic of size $2m+1$, comprised of field automorphisms. So
$\bO^{p'}(A)=A$ means $A/S$ must be a $3$-group and $B_0(A)$ is the unique block
lying above $B_0(S)$. In this case, the unique character $\chi$ of degree
$q^2(q^4-q^2+1)$ has degree divisible by $p$ and must be rational and invariant
under~$A$. Then we are done by Lemma \ref{lem:oddext}(1).

Next, let $S=\PSL_2(q)$ with $q=p^f$. Here Sylow $p$-subgroups of $A$ are
abelian unless $A$ induces field automorphisms of order $p$, thus $p|f$.
Let $r$ be a prime divisor of $(p^{2f}-1)/(p^{2f/p}-1)$, $s\in\SL_2(q)$ an
element of order~$r$ and $\chi\in\Irr(\PGL_2(q))$ the corresponding (semisimple)
Deligne--Lusztig character. By construction, $\chi$ is of $p'$-degree, but not
invariant under any non-trivial subgroup of $A$ (as $s$ is not in any proper
subfield subgroup). Thus $\chi_{A\cap\wt{S}}$ induces to an irreducible
character $\hat\vhi$ of $A$, of degree divisible by $p$. If
$r\not\equiv1\pmod p$ then $\chi$ and hence $\hat\chi$ are $\Om$-invariant.
If $p|(r-1)$, then the field automorphism in $A$ of order~$p$ fuses
$\chi$ with its images under $\Om$, and again $\chi^A$ is $\Om$-invariant.

Next, consider the case $S=\PSL_3(\eps q)$ with $q=p^f$. Here
$\wt S=\PGL_3(\eps q)$. If $p$ is odd, the unique unipotent character in
$\Irr(\wt S)\setminus\{1_{\wt S}, \St_{\wt S}\}$ has degree divisible by $p$,
lies in $B_0(\wt S)$, and is rational-valued, so we are done by
Lemma~\ref{lem:oddext}. Assume $p=2$ so $q=2^f$. Then the unipotent character
of $\wt{S}:=\PGL_3(q)$ of degree $q(q+1)$ restricts irreducibly to
$S=\PSL_3(q)$, lies in $B_0(\wt{S})$, has even degree, and extends to a
character of $A\wt{S}$ whose values lie in $\QQ(\sqrt{q})\le\QQ_8$, by
Lemmas~\ref{lem:unipextgraph} and~\ref{lem:unipprinc}, and so
is $\Om$-invariant. This character restricts to an even-degree character in
$\Irr_\Om(B_0(A))$.

Now, let $S=\PSU_3(q)$ with $q=2^f$ and $p=2$. Observe that $A/(A\cap\wt{S})$
is cyclic. If $A\le\wt{S}$ then $A=S$ and the unipotent character $\vhi$ of
degree $q(q-1)$ is as claimed. If $A\not\le\wt{S}$ then $A$ induces a
graph-field automorphism. First, suppose that $A/(A\cap\wt{S})$ has order $2$,
so is generated by $F_2^f$.  By \cite[Thm~1]{DM24}, the cuspidal
unipotent character $\vhi$ extends to a character $\hat\vhi$ of $A$ with field
of values $\QQ(\sqrt{-q})$. Note that this is again $\Om$-invariant for $p=2$.
Then we may assume that $|A/(A\cap\wt{S})|\geq 4$, and hence that $A$ contains
$\ga:=F_2^{f/2}$. Consider
$s\in\SU_3(q)=[\GU_3(q),\GU_3(q)]$ with eigenvalues $\{\ze,\ze^{-1},1\}$
where $|\ze|\neq 3$ is a prime-power divisor of $q+1$. (Recall that $q\neq 2$
since $S$ is simple.) Then the semisimple character $\chi_s$ of $\GU_3(q)$ has
odd degree, is $2$-rational since $|s|$ is odd, and is trivial on the centre.
Further, each element of $\Om$ maps $\ze$ to $\ze^{\pm 1}$, so $\chi_s$ is
$\Om$-invariant. Then we view $\chi_s\in\Irr_\Om(B_0(\wt{S}))$. Now, since
$|\ze|\mid (q+1)$, we have $\ze^{2^{f/2}}\not\in\{\ze, \ze^{-1}\}$. 
In particular, $\chi_s$ is not $\ga$-invariant.
Further, $\chi_s^\ga\neq \chi_s\alpha$ for any $\alpha\in\Irr(\wt{S}/S)$ and
$\chi_s\neq \chi_s\alpha$ for any $1\neq\alpha\in\Irr(\wt{S}/S)$, as otherwise
$s^\ga z$, respectively $s z$, would be conjugate to $s$ for some
$z\in\bZ(\GU_3(q))$ by \cite[Prop.~2.5.20]{GM20} and \cite[Prop.~7.2]{taylor},
contradicting our assumptions on $|\ze|$. 
Thus $\chi_s$ restricts to an irreducible, $2$-rational, odd-degree character
of $\Irr_\Om(B_0(S))$ that is not $A$-invariant, and we are done with this case
by Lemma~\ref{lem:oddext}(3).

Now suppose $S=\type{D}_4(q)$ with $p\leq 3$ and $q=p^f$. Let $\chi$ be one of
the unipotent characters of $\wt{S}$ listed in \cite[Thm~4.5.11(b)]{GM20}:
this character lies in the principal series, has degree divisible by $p$, and
$I_{\Aut(S)}(\chi)$ contains all field automorphisms, but does not contain the
triality graph automorphism. Then
$I_{\Aut(S)}(\chi)/\wt{S}$ is a subgroup of
$\langle\tau\rangle\times\langle F_0\rangle=C_2\times C_{f}$, where $\tau$ is a
graph automorphism of order $2$ and $F_0$ is a field automorphism of order $f$.
Note that our assumption $\bO^{p'}(A)=A$ means that $A/(A\cap\wt{S})$ is a
subgroup of $\langle\tau\rangle\times\langle F_0^m\rangle$, where $m=f_{p'}$.
By Lemmas~\ref{lem:unipextgraph} and~\ref{lem:unipprinc}, $\chi$ extends to a
rational character $\chi_1$ of $\wt{S}\langle\tau\rangle$ and to a rational
character $\chi_2$ of $\wt{S}\langle F_0^m\rangle$. Since both possible choices
of $\chi_1$ must therefore be rational, we may assume
$\chi_1\in\Irr(B_0(\wt{S}\langle\tau\rangle))$. Since $F_0^m$ has $p$-power
order, we also have $\chi_2\in\Irr(B_0(\wt{S}\langle F_0^m\rangle))$.  Then
taking the unique common extension of $\chi_1$ and $\chi_2$, we see $\chi$
extends to a rational character of $B_0(I_{A\wt{S}}(\chi))$. Then the
corresponding unipotent character $\vhi:=\chi|_S$ of $S$ extends to a rational
character $\hat\vhi$ of $B_0(I_A(\vhi))$. The induced character
$\hat\vhi^A$ then also has degree divisible by~$p$ and is rational-valued.
Further, by \cite[Cor.~6.2 and Thm~6.7]{N98}, we have
$\hat\vhi^A\in\Irr(B_0(A))$.

We may now assume that $S$ is not a Suzuki or Ree group and not of types
$\type{A}_1$, $\type{A}_2$, $\tw2\type{A}_2$, nor $\type{D}_4$ when $p\le3$. In
particular, $A/(A\cap\wt{S})$ is a $p$-group. Then the unipotent character
$\vhi\in\Irr_\Om(B_0(S))$ exhibited in the proof of Theorem~\ref{thm:ext}
restricts
irreducibly from a member of $\Irr_\Om(B_0(\wt S))$ and has degree divisible
by~$p$. Further, note $I_A(\vhi)$ contains $A\cap\wt{S}$. Then the extension
$\hat\vhi\in\Irr_\Om(B_0(I_A(\vhi)))$ from Theorem~\ref{thm:ext} has degree
divisible by $p$ and the induced character $\hat\vhi^A$ lies in
$\Irr_\Om(B_0(A))$ and has degree divisible by $p$, as required.
\end{proof}

\begin{prop}   \label{prop:as2 cross ab}
 The conclusion of Theorem~\ref{thm:as2} holds for $S$ of Lie type in
 characteristic not~$p$ when Sylow $p$-subgroups of $S$ are abelian or when
 $p=3$ and $S=\PSL_3(\eps q)$ with $q\equiv \eps\pmod 3$.
\end{prop}

\begin{proof}
In these cases, we will see that the characters with positive height in
$B_0(A)$ constructed in \cite[Props~3.10 and~3.11]{MMR} will satisfy our
conditions.
 
Let $\bG$ be of simply connected type such that $S=G/\zent{G}$ with $G=\bG^F$.
We may also identify $\wt{S}$ with $\bG_\ad^F$, with $\bG_\ad$ the
corresponding group of adjoint type. 

First, assume that $p\geq 5$. In this case, since Sylow $p$-subgroups of $S$
are abelian, we have $|\zent{G}|$ is not divisible by $p$. In
\cite[Prop.~3.11]{MMR}, a character of $\Irr(B_0(A))$ is constructed with
height 1. This is done by constructing a semisimple character
$\chi_s\in\Irr_{p'}(B_0(G))$ trivial on $\zent{G}$ and indexed by a
$p$-element $s\in G^\ast$ that is not invariant under $A$. 
Since $|s|$ is a power of $p$ and hence relatively prime to $|\zent{G}|$, we
have $\cent_{\bG^\ast}(s)$ is connected, so $\chi_s$ is the unique semisimple
character in its Lusztig series and, as an element of $\Irr(S)$, is the
restriction of a semisimple character $\wt{\chi}_{\wt s}$ of $\wt{S}$ with
$|\wt s|$ also of $p$-power order. 
Then we have $\wt\chi_{\wt s}^\sigma=\wt \chi_{\wt s}$ for any $\sigma\in\Om$,
since such $\sigma$ stabilize $p$th roots of unity.  Then we are done by
Lemma~\ref{lem:oddext}(3).

If $p=2$, then $S=\PSL_2(q)$ with $q\equiv\eps 3\pmod 8$ for $\eps\in\{\pm1\}$,
whose group of field automorphisms has odd order. Then $A=\wt{S}=\PGL_2(q)$,
and again the character of $A$ discussed in \cite[Prop.~3.11]{MMR}, labelled by
an element of order $4$, with even degree $q+\eps$, lies in $\Irr_\Om(B_0(A))$.

Now let $p=3$; then our assumption $S$ has abelian Sylow $3$-subgroups forces
$S=\PSL_2(q)$ or $\PSL_3(\eps q)$ with $\eps\in\{\pm1\}$.

Assume first that $S=\PSL_3(\eps q)$ with $q\equiv-\eps\pmod 3$ or
$S=\PSL_2(q)$. Here $S$ has cyclic Sylow $3$-subgroups. Since by assumption $A$
contains field automorphisms of $3$-power order, we have $q=q_0^3$ for $q_0$
some power of the defining characteristic, and we may consider the field
automorphism $F_{q_0}$ to lie in $A$. Let $\ze\in\FF_{q^2}^\times$ be a
3-element of order dividing $(q^2-1)_3$ but not $(q_0^2-1)_3$. Let $\wt\chi$ be
the semisimple character of $\wt{S}$ corresponding to a semisimple element
$s\in\SL_3(\eps q)$, resp.~$\SL_2(q)$ with eigenvalues $\{\ze,\ze^{-1},1\}$,
resp.~$\{\ze,\ze^{-1}\}$.  Then $\wt\chi|_S$ is irreducible and $\wt\chi$ is
not fixed by $F_{q_0}$, and hence is not stable under $A$. However, this
character is stable under $\Om$, since any $\sigma\in\Om$ fixes $3$-power roots
of unity. Since $\wt\chi$ lies in $B_0(\wt{S})$ by \cite[Cor.~3.4]{Hiss}, we
are done in this case by Lemma \ref{lem:oddext}(3). 

Now let $S=\PSL_3(\eps q)$ with $q\equiv\eps\pmod 3$. Here again we
use the same characters in \cite[Props~3.10, 3.11]{MMR}. In this case, $A/S$
is a $3$-group.  First, assume $A$ contains non-field automorphisms.  Note that
the three characters of $3'$-degree $\frac{1}{3}(q+\eps)(q^2+\eps q+1)$ of~$S$
take values in $\QQ_3$ and hence are $\Om$-invariant. Further, it is shown in
\cite[Prop.~3.10]{MMR} that $[A: I_A(\chi)]=3$ if $\chi$ is one of these
characters. Then by the same arguments in Lemma \ref{lem:oddext}(3), there is a
character of degree divisible by $3$ in $\Irr_\Om(B_0(A))$.

So now assume that $A$ contains only field automorphisms. Again the proof of
\cite[Prop.~3.10]{MMR} yields a semisimple character
$\chi:=\chi_s\in\Irr(B_0(S))$, for $s$ a 3-element, that has degree divisible
by~$3$, is invariant under $A$, and restricts irreducibly from a semisimple
character $\wt\chi$ of $\wt{S}$. Then $\chi$ and $\wt\chi$ are $\Om$-invariant
as before. By \cite[Prop.~6.7]{ruhstorfer} and its proof, $\wt\chi$ extends to
an $\Om$-invariant character of $\wt{S}A$, since $A$ is generated by $S$ and
field automorphisms.  (Indeed, note that letting $A=S\langle F'\rangle$ for a
field automorphism~$F'$ of $3$-power order, the second paragraph of the proof
of \cite[Thm.~7.4]{ruhstorfer} yields that we may apply the proof of
\cite[Prop.~6.7]{ruhstorfer} to each $\sigma\in\Om$. While our $\Om$ is not in
the group $\mathcal{H}$ there, since $\wt\chi$ is indexed by a semisimple
$3$-element, and therefore stable under $\Om$, we may still apply the arguments
there with $x_\sigma=\sigma\wt t^{-1}$ in the notation of loc.~cit.)
Then $\chi$ extends to an $\Om$-invariant character $\hat\chi$ of $A$.  Since
$A/S$ is a $3$-group, $B_0(A)$ is the unique 3-block above $B_0(S)$, and
therefore $\hat\chi$ is a character with degree divisible by $3$ in
$\Irr_\Om(B_0(A))$, as desired.
\end{proof}

\begin{prop}   \label{prop:as2 cross}
 The conclusion of Theorem~\ref{thm:as2} holds for $S$ of Lie type in
 characteristic not~$p$ when Sylow $p$-subgroups of $S$ are non-abelian.
\end{prop}

\begin{proof}
Let $S$ be of Lie type in characteristic $r\neq p$. Notice that when
$p\geq 5$, or when $p=3$ and $S$ is not of type $\type{D}_4$, we will be done
by Lemma \ref{lem:oddext}(1) if $B_0(\wt{S})$ contains a rational unipotent
character of degree divisible by~$p$ .

If $p\geq 5$, then by the first paragraph of the proof of \cite[Thm~2.18]{KM17}
there exists a unipotent character $\chi\in\Irr(B_0(\wt{S}))$ of degree
divisible
by~$p$. By \cite[Lemma~2.8]{KM17} this continues to hold when $p=3$ except when
$S=\PSL_3(\eps q)$ with $q\equiv\eps\pmod 3$. However, we may assume $S$ is not
the latter case, by Proposition~\ref{prop:as2 cross ab}. 
If $S$ is of classical type, $\chi$ is rational, see \cite[Cor.~4.4.24]{GM20}.
If $S$ is of exceptional type, by our assumption on $p$ we either have $p=3$,
or $p=5$ and $S$ is of type $\type{E}$, or $p=7$ and $S=\type{E}_8(q)$. Assume
$p=3$. Here,
$$\rho_2',\phi_{9,2},\phi_{81,6},\phi_{9,6}',\phi_{27,2},\phi_{567,6}$$
are rational unipotent characters in $B_0(S)$ of degree divisible by~3 for
$$S=\tw2F_4(q^2),F_4(q),E_6(q),\tw2E_6(q),E_7(q),E_8(q)$$
respectively, by \cite{En00} and \cite[Cor.~4.5.6]{GM20}, and
$$G_2[1],\phi_{2,1},\tw3D_4[1],\phi_{2,2}$$
are such rational unipotent characters for $G_2(q)$ with $q\equiv1,-1\pmod3$,
$\tw3D_4(q)$ with $q\equiv1,-1\pmod3$ respectively. Similarly, for the primes
$p=5,7$ an inspection of the tables in \cite[13.9]{Ca85} and \cite{En00} shows
the existence of a rational unipotent character as desired. 

If $S=\type{D}_4(q)$ and $p=3$, the same unipotent character $\chi$ of $\wt{S}$
considered above in the case of defining characteristic has degree divisible
by~$3$ in this case. Further, choosing $\chi$ more specifically with symbol
$2\choose 2$, we see this character lies in $B_0(\wt{S})$, since the $1$-core
and $1$-cocore of this symbol are both
trivial in this case. Then the same considerations as above yield a rational
character in $\Irr(B_0(A))$ of degree divisible by~$p$ above $\vhi:=\chi|_S$.
This completes the proof when $p\ge3$.

Finally, assume $p=2$. The Ree groups $\tw2\type{G}_2(q^2)$ have abelian Sylow
$2$-subgroups so do not occur here. Consider the case $S=\type{G}_2(3^n)$ with
$n\geq 1$, so that $\Aut(S)$ contains an exceptional graph automorphism $\tau$,
and suppose that $A$ contains $\tau F_1$ for some (possibly trivial) field
automorphism $F_1$. Then the character $\chi=\phi_{1,3}'$ in the notation of
\cite[Sec.~13.9]{Ca85} has odd degree, is rational-valued by
\cite[Cor.~4.5.6]{GM20}, lies in $B_0(S)$ by \cite[Thm~A]{En00}, and is fixed
by the field automorphisms but not by $\tau$ by \cite[Thm~2.5]{Ma08}. Then
$A\neq I_A(\chi)$, and by Lemma~\ref{lem:oddext}(3), our statement holds in
this case. Then we may assume that $A/S$ is comprised of field automorphisms in
the case that $S=\type{G}_2(3^n)$.

Now, assume for the moment that
$S\neq\PSL_2(q), \PSL_3(\eps q)$. Let $\chi$ be a rational-valued, principal
series unipotent character in $B_0(S)$ with even degree guaranteed by
Proposition~\ref{prop:evenrat}.  In the case that $S=\PSL_n(\eps q)$ with
$n\geq 4$, let $\la$ be the partition $(n-2, 2)$ if $n\equiv 0,3\pmod 4$ and
$\la=(n-2, 1^2)$ if $n\equiv 1,2\pmod 4$. Then taking $\chi$ more specifically
to be the character indexed by $\la$, we have $\chi$ is such a character using
the degree formula in \cite[Prop.~4.3.2]{GM20}, but also does not satisfy the
condition in Lemma~\ref{lem:unipextgraph}. Then in all relevant cases, $\chi$
extends to a rational-valued character of $\wt{S}\langle\tau\rangle$ if $\tau$
is a non-trivial graph automorphism stabilizing~$\chi$, by
Lemma~\ref{lem:unipextgraph}.

If $A/(A\cap \wt{S})$ is abelian, note that it must be a $2$-group, using our
assumption that $\bO^{2'}(A)=A$. Then arguing as before, with
Lemmas~\ref{lem:unipextgraph} and~\ref{lem:unipprinc} and unique common
extensions, we obtain an extension $\hat\vhi$ of $\vhi:=\chi|_S$ in
$\Irr_\Om(B_0(I_A(\vhi)))$. Again applying \cite[Cor.~6.2 and Thm~6.7]{N98}, we
have $\hat\vhi^A\in\Irr_\Om(B_0(A))$ with even degree.
 
Next, we assume that $A/(A\cap\wt{S})$ is non-abelian, so $S=\type{D}_4(q)$.
We choose $\chi$ in this case to be the unique (unipotent) character of degree
$q(q^2+1)^2$, so that $\chi$ extends to a rational character of~$\Aut(S)$ by
the second paragraph of \cite[Sec.~4.5]{TT22}, To find a character in $B_0(A)$,
note that $A/(A\cap\wt{S})\leq X:= \sym_3\times\langle F_0^m\rangle$, where
$F_0$ is a field automorphism of order $f$ and $f_{2'}=m$, for $q=r^f$. Since
every
character of $X$ is $\Om$-invariant (as the characters of $\sym_3$ are rational
and $\Om$ stabilizes $2$-power roots of unity), it follows by Gallagher's
theorem that every character of $A$ lying above $\chi':=\chi_{A\cap\wt{S}}$ is
$\Om$-invariant. In particular, we know that there is a character above $\chi'$
in $B_0(A)$, which therefore has even degree and lies in $\Irr_\Om(B_0(A))$.
 
It remains to discuss the groups $\PSL_2(q)$ and $\PSL_3(\eps q)$. First, let
$S=\PSL_2(q)$. Since Sylow 2-subgroups of $S$ are non-abelian by assumption, we
have $q\equiv\pm1\pmod8$. Here
$A/(A\cap\wt{S})$ is generated by field automorphisms of 2-power order. Let
$s\in\SL_2(q)$ be a 2-element of maximal 2-power order. Then the corresponding
Deligne--Lusztig character $\chi$ of $\wt{S}=\PGL_2(q)$ of degree $q\pm1$ has
values in $\QQ_{2^f}$ for some $f\ge1$, is not invariant under any field
automorphism of 2-power order, and it restricts irreducibly to $S$. Thus
its irreducible induction to $A$ is as required. 

Finally, let $S=\PSL_3(\eps q)$ with $q\not\equiv\eps\pmod 3$ odd, and continue
to assume $p=2$. In the first paragraph of \cite[Prop.~2.12]{MN}, a character
$\psi\in\Irr_\Om(B_0(S))$ is constructed such that $I_{\Aut(S)}(\psi)=\wt{S}$
and $\psi=\chi_S$ for a semisimple $\chi\in\Irr_\Om(B_0(\wt{S}))$ with even
degree corresponding to a suitable $2$-element. (We note that, although only
a specific $\sigma\in\Om$ is considered there, $\chi$ is $\Om$-stable for the
same reason --- that any such Galois automorphism fixes $2$-power roots of
unity). Then as in loc.~cit., the induced character $(\chi_{\wt{S}\cap A})^A$
satisfies our statement.
\end{proof}

%%%%%%%%%%%%%%%%%%%%%%%%%%%%%%%%%%%%%%%%%%%%%%%%%%%%%%%%%%%%%%%%%%%%%%%%%
\section{The Galois It\^o--Michler theorem}   \label{sec:IM}

In this section we prove Theorem \ref{im}.
We will use the Alperin--Dade character correspondence.

\begin{lem}   \label{lem:alperin-dade}
 Suppose that $N$ is a normal subgroup of $G$, with $G/N$ a $p'$-group.
 Let $P \in \Syl_p(G)$ and assume that $G=N\bC_G(P)$. Then restriction of
 characters defines a natural bijection between the irreducible characters of
 the principal $p$-blocks of $G$ and $N$. In particular restriction defines a
 bijection
 $$\res: \Irr_\Om(B_0(G))\rightarrow\Irr_\Om(B_0(N)).$$
\end{lem}

\begin{proof}
The case where $G/N$ is solvable was proved in \cite{Alp} and the general case
in \cite{Dad}. The last part of the statement follows immediately since $\tau$
acts on $\Irr(B_0)$ for every $\tau\in\cG$.
\end{proof}

If $N\trianglelefteq G$ and $\theta\in\Irr(N)$, we denote by $\Irr(G|\theta)$
the set of irreducible characters of $G$ lying over $\theta$, and by
$\cd(G|\theta)$ the set of degrees of the characters in $\Irr(G|\theta)$.

\begin{lem}   \label{up}
 Suppose that $N$ is a normal subgroup of $G$, with $G/N$ a $p'$-group. Let
 $\Om=\Om(|G|)$. If $\theta\in\Irr_\Om(N)$ then there exists
 $\chi\in\Irr_\Om(G)$ over $\theta$. Furthermore, if 
 $\theta\in\Irr_\Om(B_0(N))$ then there exists $\chi\in\Irr_\Om(B_0(G))$ over
 $\theta$. 
\end{lem}

\begin{proof}
We argue by induction on $|G:N|$. Suppose that $G_{\theta}<G$. By the inductive
hypothesis there exists $\psi\in\Irr_\Om(G_\theta|\theta)$. By Clifford's
correspondence, $\chi=\psi^G$ is irreducible and by the formula for the induced
character $\QQ(\chi)\subseteq\QQ(\psi)$ so $\chi$ is $\Om$-invariant. Moreover,
if $\psi\in\Irr(B_0(G_\theta))$, then $\chi\in\Irr(B_0(G))$ by \cite[Cor.~6.2
and Thm~6.7]{N98}. Hence, we may assume that $\theta$ is $G$-invariant in both
statements. 

For $d\ge1$ set $\Irr_d(G|\theta):=\{\chi\in\Irr(G|\theta)\mid \chi(1)=d\}$.
Since $\theta$ is $\Om$-invariant, $\Om$ acts on $\Irr_d(G|\theta)$ for every
$d$. Since
$$|G:N|=\sum_{\chi\in\Irr(G|\theta)}\left(\frac{\chi(1)}{\theta(1)}\right)^2
  =\sum_{d}|\Irr_d(G|\theta)|\frac{d^2}{\theta(1)^2}$$
is a $p'$-number, we conclude that there exists $d\in\cd(G|\theta)$ such that
$|\Irr_d(G|\theta)|$ is not divisible by $p$. Since $|\Om|$ is a power of $p$,
it follows that there exists $\chi\in\Irr_d(G|\theta)$ that is $\Om$-invariant.

Now assume that $\theta\in\Irr_\Om(B_0(N))$. If $G/N$ is not simple, let
$M\trianglelefteq G$ with $N<M<G$. By induction there is
$\psi\in\Irr_\Om(B_0(M))$ lying over $\theta$. Again by induction, there is
$\chi\in\Irr_\Om(B_0(G))$ lying over $\psi$ (and hence over $\theta$). Hence we
may assume that $G/N$ is simple.

Now, let $P\in\Syl_p(G)$ and notice that $P\leq N$. By the Frattini argument,
$G=N\bN_G(P)$ and hence $M:=N\bC_G(P)\trianglelefteq G$. Since $G/N$ is simple
we have that $M=N$ or $M=G$. If $M=G$, by Lemma~\ref{lem:alperin-dade} there
exists $\chi\in\Irr_\Om(B_0(G))$ such that $\chi_N=\theta$, and we are done.
Hence we may assume that $M=N$. In this case, $\bC_G(P)\leq N$ and by
\cite[Lemma~4.2]{MMR} we have that $B_0(G)$ is the only block covering
$B_0(N)$. By the first part of the proof, there is $\chi\in\Irr_\Om(G|\theta)$,
hence $\chi\in\Irr_\Om(B_0(G)|\theta)$ and we are done. 
\end{proof}

The following is the $p$-group case of Theorem~\ref{im}.

\begin{lem}   \label{lem:itomichlerpgroup}
 Let $p$ be a prime number and let $P$ be a $p$-group. Then
 $$\Irr(P)=\Irr_\Om(P).$$
\end{lem}

\begin{proof}
By the definition of $\Om$, every $\tau\in\Om$ fixes the $p$-power roots of
unity, so the result is clear.
\end{proof}

We need the following technical lemma.

\begin{lem}   \label{lem:pn}
 Let $p$ and $q$ be distinct primes. Suppose that $G=VP$, where $P$ is an
 abelian $p$-group, $V$ is an irreducible $P$-module over $\FF_q$ and
 $\bC_P(V)=1$. Suppose that $p$ divides $q-1$. Then there exists
 $\chi\in\Irr(G)$ of degree divisible by $p$ such that
 $\QQ(\chi)\subseteq\QQ_q$ and $p$ divides $|\QQ_q:\QQ(\chi)|$.
\end{lem}

\begin{proof} 
Since $V$ is a faithful irreducible $P$-module, and $\bC_V(x)\trianglelefteq G$
for every $x\in P$, we have that $\bC_G(x)=1$ for every $x\in P$. 
Thus $G$ is a Frobenius group with abelian complement $P$.  Let $h\in P$ be an
element of order $p$ and put $H=\langle h\rangle$. 

Now, write $V=V_1\oplus\cdots\oplus V_t$ as a direct sum of irreducible
$H$-modules. Since $p$ divides $q-1$, the field $\FF_q$ contains primitive
$p$th roots of unity, so all irreducible $H$-modules have dimension $1$.
Since $\bC_P(V)=1$, there exists $i$ and $\la\in\Irr(V_i)$ that is not
$H$-invariant by \cite[Prop.~12.1]{MW}. Write $V=U\oplus W$, so that $U=V_i$
and $W\trianglelefteq VH$.  We view $\la$ as a character of $V/W$ 
(hence of $V$) that induces irreducibly to $\mu=\la^{VH}\in\Irr(VH)$.  Let
$U=\langle v\rangle$, so that $\la(v)=\veps$ for some primitive $q$th root of
unity $\veps$ (note that since $\la$ is not $H$-invariant $\la(v)\neq1$).
Since $H$ does not act trivially on $U$, $v^{h^{-1}}=v^i$ for some~$i$ that has
order~$p$ modulo~$q$.  Let $\tau\in\Gal(\QQ_q/\QQ)$ be the Galois automorphism
such that $\tau(\veps)=\veps^i$, so that $\la^h(v)=\la(v)^{\tau}$. Note that
$\tau$ has order $p$ since $\tau^p(\veps)=\veps^{i^p}=\veps$. Note that
$$\mu_V=\sum_{j=0}^{p-1}\la^{g^j}$$
so that
$$\mu(v)=\veps+\veps^i+\cdots+\veps^{i^{p-1}}
  =\veps+\veps^\tau+\cdots+\veps^{\tau^{p-1}}=\mu^\tau(v).$$
Thus $\mu$ is $\tau$-invariant and hence $\QQ(\mu)$ is $\tau$-invariant, that
is, $\tau\in\Gal(\QQ_q/\QQ(\mu))$. Since $\tau$ has order $p$, the fundamental
theorem of Galois theory implies that $|\QQ_q:\QQ(\mu)|$ is divisible by $p$.

Since $G$ is a Frobenius group, $\chi=\la^G\in\Irr(G)$ (so $p$ divides
$\chi(1)$). Note that $\QQ(\chi)=\QQ(\chi_V)$ (because $\chi$ vanishes
off~$V$). By Clifford's theorem,  $\chi_{VH}$ is a sum of conjugates of $\mu$.
Therefore,
$$\QQ(\chi)=\QQ(\chi_{VH})\subseteq\QQ(\mu).$$ 
It follows that $|\QQ_q:\QQ(\chi)|$ is divisible by $p$, as wanted.
\end{proof}

Given a group $G$, we write $\bO_p(G)$ to denote the largest normal
$p$-subgroup of $G$. Analogously $\bO_{p'}(G)$ is the largest normal
$p'$-subgroup of $G$. The following is Theorem~\ref{im}:

\begin{thm}   \label{thm:im}
 Suppose that $G$ does not have composition factors isomorphic to $S$ with
 $(S,p)$ one of the pairs listed in Theorem~\ref{thm:as}. Then all characters
 in $\Irr_\Om(G)$ have $p'$-degree if and only if $G$ has a normal abelian
 Sylow $p$-subgroup.
\end{thm}

\begin{proof}
The ``if'' direction is clear, so we prove the ``only if'' direction. Let $P$
be a Sylow $p$-subgroup of $G$. 
\medskip

{\it Step 1.} If $P$ is normal in $G$, then $P$ is abelian. 
\medskip

Let $\theta\in\Irr(P)$. By Lemma~\ref{lem:itomichlerpgroup} we have that
$\theta\in\Irr_\Om(P)$. By Lemma~\ref{up}, there exists $\chi\in\Irr_\Om(G)$
over $\theta$. By hypothesis, $\chi$ has $p'$-degree. Since $\chi(1)$ is a
multiple of $\theta(1)$, this implies that $\theta(1)$ is $p'$ and hence
$\theta$ is linear. Thus $P$ is abelian.
\bigskip

In the following we want to prove that $P$ is normal in $G$. Let $G$ be a
minimal counterexample. 
\bigskip

{\it Step 2.} If $1<N\trianglelefteq G$, $G/N$ has normal and abelian Sylow
$p$-subgroups. In particular, there is a unique minimal normal subgroup of $G$
and $\bO_{p}(G)=1$.
\medskip

Since $\Irr_\Om(G/N)\subseteq\Irr_\Om(G)$, the hypothesis is inherited by
$G/N$. By the minimality of $G$ as a counterexample, $G/N$ has a normal and
abelian Sylow $p$-subgroup. Now, if
$M,N$ are minimal normal subgroups, $G$ is isomorphic to a subgroup of
$G/N\times G/M$, and hence $G$ has a normal and abelian Sylow subgroup, as
wanted. Finally, if $1<N=\bO_p(G)$, $P/N$ is normal in $G$, so $P$ is normal
in $G$, and we are done. 
\medskip

{\it Step 3:} We have that $G=\bO^{p'}(G)$. In particular, if $N$ is the
unique minimal normal subgroup of $G$, then $G=PN$.
\medskip

Let $L=\bO^{p'}(G)$. Suppose that $L<G$.  Let $\theta\in\Irr_\Om(L)$. By
Lemma~\ref{up}, there exists $\chi\in\Irr_\Om(G)$ over $\theta$. By hypothesis,
$\chi(1)$ is $p'$. We conclude that all characters in $\Irr_\Om(L)$ have
$p'$-degree. By induction, $L$ has a normal Sylow $p$-subgroup and hence the
same holds for $G$. Thus, we may assume that $G=\bO^{p'}(G)$. By Step 2, $PN/N$
is normal in $G/N$, so $PN$ is normal in $G$, which forces $G=PN$.
\medskip

{\it Step 4:} Let $N$ be the unique minimal normal subgroup of $G$. Then $N$ is
a direct product of non-abelian simple groups.
\medskip

Suppose that $N$ is elementary abelian. By Step 2, $N$ is a $q$-group for some
prime $q\neq p$ and by Step~3, $G=PN$. Notice that $P\cong G/N$ is abelian by
Step~2. Then $\bC_P(N)$ is normal in $G$, so $\bC_P(N)=1$ again by Step~2. 
Note that $\Om(|G|)$ is cyclic.
If $p$ does not divide $q-1$ then $\Om$ is trivial and we finish the proof
using the It\^o--Michler theorem. Hence, we may assume that $p$ divides $q-1$,
so that $\Om$ is cyclic of order $p$.
By Lemma \ref{lem:pn}, there exists $\chi\in\Irr(G)$ of degree divisible by $p$
and such that $p$ divides $|\QQ_q:\QQ(\chi)|$. Then there is
$\tau\in\Gal(\QQ_q/\QQ(\chi))$ of order $p$. It follows that
$\Om$ is generated by an extension of $\tau$, whence $\chi\in\Irr_\Om(G)$ has
degree divisible by~$p$, a contradiction. 
\medskip

{\it Step 5:} Completion of the proof. 
\medskip

Let $N$ be the unique minimal normal subgroup of $G$. By Step 4, we know that
$N$ is a direct product of isomorphic non-abelian simple groups. Write
$N=S_1\times\cdots\times S_t$, where $S_i\cong S$ is non-abelian simple. By
Theorem~\ref{thm:as}, we may assume that $t>1$. Write
$H=\bigcap_{i=1}^t\bN_G(S_i)$
so that $G/H$ is isomorphic to a transitive permutation group on $t$ letters
and $N\leq H\leq \Aut(S)\times\cdots\times\Aut(S)$. By Step~2 and Step~3,
$G/H$ is an abelian $p$-group. Therefore, all point stabilizers are trivial. 

By Theorem~\ref{thm:ext}, there exists $1_S\neq\vhi\in\Irr(S)$ that extends to
an $\Om$-invariant character of its inertia group in $\Aut(S)$. Since
$G_\nu\leq H$, $\nu=\vhi\times1_S\times\cdots\times1_S$ extends to an
$\Om$-invariant character $\tilde{\nu}$ of $G_\nu<G$. Hence, since $G/H\ne1$,
$\tilde{\nu}^G\in\Irr(G)$ has degree divisible by $p$ and is $\Om$-invariant.
This contradicts the hypothesis. 
\end{proof}

The following is the $p$-solvable case of Theorem \ref{main}. 

\begin{cor}   \label{psol}
 Let $p$ be a prime. Let $G$ be $p$-solvable and let $P\in\Syl_p(G)$.
 Then all characters in $\Irr_\Om(B_0(G))$ have height zero if and only if $P$
 is abelian.
\end{cor}

\begin{proof}
Since $G$ is $p$-solvable, $\Irr(B_0(G))=\Irr(G/\bO_{p'}(G))$. Now, the result
follows from Theorem~\ref{thm:im} and Hall--Higman's Lemma 1.2.3.
\end{proof}

\begin{rem}   \label{alt}
We remark that if $p$ is odd, our arguments can be adapted to show that we can 
replace $\Om$ by $\cP=\{\sigma\in\Gal(\QQ_{|G|}/\QQ)\mid o(\sigma)=p\}$ in
Theorem~\ref{thm:im} and Corollary~\ref{psol} (using that the field of values
of an irreducible character of an odd order $p$-group is a full cyclotomic
field). This is not possible however in Theorem~\ref{main}. 
Fix $p=3$ and let $G=X\star Y$ be the central product of a cyclic group $X$ of
order $9$ and the $3$-fold cover $Y$ of the alternating group on $6$ letters.
Note that $G$ does not have abelian Sylow $3$-subgroups. We claim that every
character in $\Irr_{\cP}(B_0(G))$ has $3'$-degree. We can check in \cite{Atl}
that every character of degree divisible by $3$ of $B_0(Y)$ has field of values
of $2$-power degree. Hence, all of these characters are invariant under a Sylow
$3$-subgroup of the Galois group. Hence, $B_0(G)=B_0(X)\star B_0(Y)$ (see
\cite[Lemma~4.1]{MN}) does not have $\cP$-invariant characters of degree
divisible by $3$, as claimed.
\end{rem}

We conclude this section with a theorem that generalizes results of
Dolfi--Navarro--Tiep \cite[Thm~A]{DNT}, who consider the case that $\sigma$ is
complex conjugation, as well as the $p=2$ case of Grittini \cite[Thm~A]{gri},
which considers the same statement for $p$-solvable groups (there for an
arbitrary prime $p$). This is also part of \cite[Thm~A]{Gr}.

\begin{thm}   \label{thm:evenorder2}
 Let $\sigma\in\Gal(\QQ_n/\QQ)$ have order~$2$ and let $G$ be a finite group of
 order dividing~$n$. If every $\chi\in\Irr(G)$ fixed by $\sigma$ has odd
 degree, then $G$ has a normal Sylow $2$-subgroup.
\end{thm}

\begin{proof}
We proceed by induction on $|G|$.  Using \cite[Prop.~2.3]{gri} and following
Steps 1--2 of the proof of \cite[Thm~A]{gri}, we see that we may assume
$G=\bO^{2'}(G)=NP$, where $N$ is the unique minimal normal subgroup of $G$ and
$P$ is a Sylow $2$-subgroup of $G$. 

Suppose first that $N$ is abelian. Then the result follows from
\cite[Thm~A]{gri}.

Hence, we assume that $N=S^{x_1}\times\cdots\times S^{x_t}$ is the product of
conjugates of a non-abelian simple subgroup $S$ and let
$H:=\norm_G(S)/\cent_G(S)$. Now, by \cite[Thm~C]{TT22} and
Theorem~\ref{thm:as}, there is a $\sigma$-invariant character $\chi\in\Irr(H)$
of even degree, and such that there is some non-trivial $\theta\in\Irr(S)$
lying under $\chi$. From here, we argue exactly as in the proof of
\cite[Thm~A]{DNT}. Namely, taking
$\eta:=\theta\times 1_{S^{x_2}}\times\cdots\times 1_{S^{x_t}}\in\Irr(N)$ we have
$I_G(\eta)\le H$.  Since $\chi$ (now viewed as the inflation to $\norm_G(S)$)
lies over $\eta$, we obtain that the induced character $\chi^G$ is a
$\sigma$-invariant, irreducible character of even degree, a contradiction. 
\end{proof}

We care to remark that the proofs in this section (with minor modifications in
Lemma~\ref{lem:pn}), as well as the ones in Section~4, would still work if we
replace $\cJ$ by the bigger group of Galois automorphisms of $p$-power order
that fix $p$-power roots of unity mentioned in the introduction.

%%%%%%%%%%%%%%%%%%%%%%%%%%%%%%%%%%%%%%%%%%%%%%%%%%%%%%%%%%%%%%%%%%%%%%%%%
\section{The Galois height zero conjecture}   \label{sec:bhz}

We assume that the reader is familiar with the basic properties of the
generalized Fitting subgroup $\bF^*(G)$ (see e.g. \cite[Sec.~9A]{isagr}). We
write $\bE(G)$ to denote the layer of $G$ and $\bF(G)$ to denote the Fitting
subgroup of $G$. 

Now, we prove the main result, which we restate:

\begin{thm}
 Let $G$ be a finite group, let $p$ be a prime number and let $P\in\Syl_p(G)$.
 Suppose that $G$ does not have composition factors isomorphic to $S$ with
 $(S,p)$ one of the pairs listed in Theorem \ref{thm:as}. Then all characters
 in $\Irr_\Om(B_0(G))$ have height zero if and only if $P$ is abelian.
\end{thm}

\begin{proof}
The ``if" part follows from the ``if" part of Brauer's height zero conjecture,
proved in \cite{KM}. We prove the ``only if" part. 
Let $P$ be a Sylow $p$-subgroup of $G$. We want to see that $P$ is abelian.
Let $G$ be a minimal counterexample. 
\medskip

{\it Step 1:} We have that $\bO_{p'}(G)=1$. Furthermore, for every
$1<N\trianglelefteq G$, $G/N$ has abelian Sylow $p$-subgroups and $G$ has a
unique minimal normal subgroup.
\medskip

Put $M=\bO_{p'}(G)$ and assume that $M>1$.
Since $\Irr_\Om(B_0(G/M))\subseteq\Irr_\Om(B_0(G))$, all the characters in
$\Irr_\Om(B_0(G/M))$ have $p'$-degree. By the minimality of $G$ as a
counterexample, $G/M$ and hence $G$, have abelian Sylow $p$-subgroups. The
second part is analogous. For the last part, suppose that $N,M$ are minimal
normal subgroups of~$G$. Then $G$ is isomorphic to a subgroup of
$G/N\times G/M$ and the conclusion holds. 
\medskip

{\it Step 2:} We have that $G=\bO^{p'}(G)$. 
\medskip

Write $M=\bO^{p'}(G)$ and assume that $M<G$. Let $\psi\in\Irr_\Om(B_0(M))$.
By Lemma~\ref{up}, there exists $\chi\in\Irr_\Om(B_0(G))$ over $\psi$ of
$p'$-degree. Hence $\psi$ has $p'$-degree. By the minimality hypothesis $M$ has
abelian Sylow $p$-subgroups and hence the same holds for $G$. 
\bigskip

In the following, let $N$ be the unique minimal normal subgroup of $G$.
Set $K/N=\bO_{p'}(G/N)$, so that by Step~1, Step~2 and \cite[Thm~4.1]{MMR},
$G/K=X/K\times Y/K$, where $X/K$ is an abelian $p$-group and $Y/K$ is a direct
product of non-abelian simple groups of order divisible by $p$ with abelian
Sylow $p$-subgroups.

We first consider the case when $N$ is an elementary abelian $p$-group.
If $Y=K$, $G$ is $p$-solvable, so it has just one $p$-block by Step~1. In this
case the result holds by Corollary~\ref{psol}. So we may assume that $Y>K$.
\bigskip

{\it Step 3:} Let $M=\bC_G(N)$. Then $G/M$ is an abelian $p$-group.
\medskip

Let $Q$ be a Sylow $p$-subgroup of $M$, so that $\bC_G(Q)\subseteq\bC_G(N)=M$
(note that since $N$ is a normal $p$-subgroup of $M$, $N\subseteq Q$).
By \cite[Lemma~4.2]{MMR}, $\Irr(G/M)\subseteq\Irr(B_0(G))$, so 
$\Irr_\Om(G/M)\subseteq\Irr_\Om(B_0(G))$. By hypothesis, this implies that all
characters in $\Irr_\Om(G/M)$ have $p'$-degree, then by Theorem~\ref{thm:im}
and Step~2 we obtain that $G/M=PM/M$ is abelian. 
\medskip

{\it Step 4:} We have that $K>N$.
\medskip

Suppose that $K=N$, so $G/N=X/N\times Y/N=X/N\times S_1/N\times \cdots\times S_t/N$, where $S_i/N$ is non-abelian simple of order divisible by $p$ for every $i=1,\dots,t$. Now, $X\ge N$ is a $p$-group, so $1<\bZ(X)\lhd G$ and then $N\subseteq \bZ(X)$. Hence $X\subseteq \bC_G(N)=M$. Since $G/X\cong Y/N$ does not have normal subgroups of $p$-power index, we necessarily have that $M=G$. Then $N\subseteq \bZ(G)$ and $|N|=p$. Now we have that $G$ is the central product of $X$ and $S_1,\ldots, S_t$, where $Y/N\cong S_1/N\times\cdots\times S_t/N$. Since $N$ is the unique minimal normal subgroup of $G$, $N\subseteq S_i'$. Hence 
all $S_i$ are perfect, so $S_i$ is quasi-simple with centre $N$, for every $i$.
Let $1_N\neq\la\in\Irr(N)$. By Theorem~\ref{thm:quasisimples} there exists $\psi_i\in\Irr_{\Om(|S_i|)}(B_0(S_i))=\Irr_\Om(B_0(S_i))$ lying over $\la$. (If necessary, replacing $\psi_i$ by a Galois conjugate.)  Now, let $\xi\in\Irr(X|\la)$, then $\xi\in\Irr_\Om(B_0(X)|\la)$. By \cite[Lemma 4.1]{MN}
the central product of characters 
$$\chi=\xi\star \psi_1\star\psi_2\star \cdots\star \psi_t$$
lies in the principal block of $G$. Hence $\chi\in\Irr_\Om(B_0(G))$ has degree
divisible by $p$, which contradicts the assumption $K=N$. 
\medskip

{\it Step 5:} We have that $\bF(G)=\bF^*(G)$.
\medskip

In this step we use arguments from the proof of Theorem~4.6 of \cite{MMR}.
By Step 1 and the assumption that $N$ is abelian, $F=\bF(G)=\bO_p(G)>1$. Suppose that $E=\bE(G)>1$ and let $Z=\bZ(E)$. Since $N$ is the unique minimal normal subgroup of $G$, $N\subseteq Z$ (notice that $Z>1$ since otherwise $\bF^*(G)=\bF(G)\times E$ in contradiction to Step 1). We claim that $E/Z=S_1/Z\times\cdots\times S_n/Z$, where $S_i\trianglelefteq G$ for every $i$. Let $W/Z$ be a non-abelian chief factor of $G/Z$ contained in $E/Z$. By the Schur--Zassenhaus theorem and Step 1, we know that $|W/Z|$ is divisible by $p$. Now, by \cite[Thm 4.1]{MMR} applied to
$G/Z$, we have that $W/Z$ is simple and the claim follows.
Write $S=S_1$, so that $S'$ is a quasi-simple normal subgroup of $G$. Using
again that $N$ is the unique minimal normal subgroup of $G$, we have that
$N\subseteq \bZ(S)\cap S'\subseteq \bZ(S')$.
Looking at the Schur multipliers of the simple groups \cite{Atl}, if $p\geq 5$,
we deduce that $\bZ(S')$ has cyclic Sylow $p$-subgroups. Arguing as in Step~4
of the proof of Theorem 4.6 of \cite{MMR} we have that
$\bZ(S')$ has cyclic Sylow $p$-subgroups for $p=2,3$ as well.
\medskip

In all cases, we conclude that $N$ is cyclic and hence $|N|=p$. Now, the order
of $G/\bC_G(N)$ divides $p-1$. By Step 2, $G=\bC_G(N)$, so $N$ is central
in~$G$. Thus $K$ is the direct product of $N$ and a $p$-complement $H$. Since
$H$ is normal in $G$, we have a contradiction with Step 1. Hence $E=1$ and
$\bF(G)=\bF^*(G)$ as wanted.
\medskip

{\it Step 6:} We have that $N=\bF(G)=\bF^*(G)$. 
\medskip

Let $H$ be a $p$-complement of $K$, so $K=HN$ and $H\cap N=1$. Then by the
Frattini argument and the Schur--Zassenhaus theorem, we have $G=N\bN_G(H)$.
Write $L=\bN_G(H)$. Now, since $N$ is abelian and normal in $G$, $\bN_N(H)$ is
normal in $G=N\bN_G(H)=NL$, and hence $\bN_N(H) =1$ or $\bN_N(H) =N$. If
$\bN_N(H)=N$, then $H\lhd G$, and we get a contradiction since $\bO_{p'}(G)=1$
by Step~1. Thus, $L\cap N=\bN_N(H) =1$ and $L$ is a complement of $N$ in $G$.

Let $F=\bF(G)=\bO_p(G)$. We will show that $F=N$. Notice that since $\bZ(F)>1$,
we have $N\subseteq\bZ(F)$. Let $F_1=F\cap L$. Then $F_1\lhd L$ and since
$G=NL$, we have $F_1\lhd G$. Since $N$ is the unique minimal normal
subgroup, this forces $F_1=1$, so $F=N$ as wanted. Since $\bF(G)=\bF^*(G)$,
the claim follows, and the proof of the Step is completed. 
\medskip

{\it Step 7:} Completion of the proof when $N$ is elementary abelian.
\medskip

Step 6 implies that $\bC_G(N)=N$. By Step 3, $G/N$ is a $p$-group, but then $G$
is a $p$-group and we are done since $\Irr_\Om(B_0(G))=\Irr(G)$ in this case
(see Lemma \ref{lem:itomichlerpgroup}).
\bigskip

Hence, from now on we may assume that $N$ is a direct product of non-abelian
simple groups.
Write $N=S_1\times\cdots\times S_t$, where $S_i\cong S$ is non-abelian simple.
By Theorem~\ref{thm:as2}, we may assume $t>1$. Write $H=\bigcap_{i=1}^t\bN_G(S)$ so
that $G/H$ is isomorphic to a transitive permutation group on $t$ letters and
$N\leq H\leq \Aut(S)\times\cdots\times\Aut(S)$
\medskip

{\it Step 8:} We have that $\bC_G(P)\subseteq H$. In particular, $G/H$ is a
non-trivial $p$-group and $G/N$ has a normal $p$-complement $K/N$. 
\medskip

Let $R\in\Syl_p(S)$, so that $Q=R\times\cdots\times R\subseteq P$ is a Sylow
$p$-subgroup of $N$.  Since $\bO_{p'}(G)=1$, $R>1$. Let $g\in G-H$. Since $g$
permutes the copies of $S$, we may assume without loss of generality that $g$
does not centralize $(x,1,\dots,1)$, where $1\neq x\in R$. The first part
follows.
\medskip

Now, using \cite[Lemma 4.2]{MMR} we have
$\Irr_\Om(G/H)\subseteq\Irr_\Om(B_0(G))$. By hypothesis, $p$ does not divide the degree of any character in $\Irr_\Om(G/H)$. It follows from Theorem~\ref{thm:im} that $G/H$ has a normal abelian Sylow $p$-subgroup. Since $\bO^{p'}(G)=G$, we conclude that $G/H$ is a (non-trivial) $p$-group, as desired.
\medskip

Finally, we prove the third claim. Since $H/N$ is isomorphic to a subgroup of $\Out(S)^t$, it follows from Schreier's conjecture that $H/N$ is solvable. Since $G/H$ is also solvable, we conclude that $G/N$ is solvable. By Step~1, $G/N$ has abelian Sylow $p$-subgroups. Therefore, by Hall--Higman's Lemma~1.2.3 applied to $G/K$ and the fact that $\bO^{p'}(G)=G$, we conclude that $G/N$ has a normal $p$-complement $K/N$.
\medskip

{\it Step 9:} We have that $G=NP$.
\medskip

In this step we follow the arguments from \cite{MN}.
Let $Q=P\cap N\in\Syl_p(N)$. By the Frattini argument, $G=N\bN_G(Q)$. Therefore, $M=N\bC_G(Q)\trianglelefteq G$. By \cite[Lemma 4.2]{MMR}, all irreducible characters of $G/M$ belong to $B_0(G)$. By Theorem~\ref{im}, $G/M$ has a normal Sylow $p$-subgroup. Since $\bO^{p'}(G)=G$, we conclude that $G/M$ is a $p$-group. Hence $K\subseteq M=N\bC_G(Q)$. Therefore, $K=N\bC_K(Q)$.
\medskip

Let $\alpha\in\Irr_\Om(B_0(NP))$. We want to show that $\alpha$ has $p'$-degree. Let $\theta\in\Irr(N)$ under $\alpha$, so $\theta$ belongs to the principal $p$-block of $N$. By Alperin's isomorphic blocks
(Lemma~\ref{lem:alperin-dade}), there exists a unique extension $\eta$ of $\theta$ in $B_0(K)$. Let $I=G_{\eta}$, so that $J=I\cap PN=(PN)_{\theta}$ (using uniqueness in Alperin's theorem). Let $\mu\in\Irr(J)$ be the Clifford correspondent of $\alpha$ over $\theta$. By the Isaacs restriction correspondence (\cite[Lemma 6.8]{N18}), let $\rho\in\Irr(I|\eta)$ be such that $\rho_J=\mu$. By the Clifford correspondence, $\chi=\rho^G\in\Irr(G)$ lies in the principal $p$-block of $G$ (because $\eta$ does and $G/K$ is a $p$-group). Now, let $\sigma\in \Om$. Since $\alpha$ is $\sigma$-invariant, we have $\theta^\sigma=\theta^{g_\sigma}$ for some $g_\sigma\in P$ by Clifford's theorem. By the uniqueness in the Alperin--Dade correspondence we have that, $\eta^{g_\sigma}=\eta^\sigma$. By uniqueness in the Clifford correspondence and the Isaacs correspondence we obtain that $\mu^\sigma=\mu^{g_\sigma}$, $\rho^\sigma=\rho^{g_\sigma}$ and $\chi^\sigma=\chi$. Since this occurs for every $\sigma\in\Om$, we conclude that $\chi\in\Irr_\Om(B_0(G))$. By hypothesis, $\chi$ has $p'$-degree. Thus $I=G$, $\chi=\rho$ and $\chi_{PN}=\alpha$ has $p'$-degree as wanted. By the minimality of $G$ as a counterexample, we may assume that $G=NP$, as claimed.
\medskip

{\it Step 10:} Completion of the proof.
\medskip

By Theorem~\ref{thm:ext}, there exists $1_S\neq\vhi\in\Irr_\Om(B_0(S))$ that
extends to an $\Om$-invariant character of the principal $p$-block of its
inertia group in $\Aut(S)$. Therefore,
$\nu=\vhi\times1_S\times\cdots\times1_S\in\Irr(B_0(N))$ extends to an
$\Om$-invariant character $\tilde{\nu}$ of the principal $p$-block of
$G_\nu\subseteq H$ (using again that $G/H$ is an abelian transitive permutation
group). Hence, $\tilde{\nu}^G\in\Irr(B_0(G))$ has degree divisible by $p$ and
is $\Om$-invariant. This contradicts the hypothesis.
\end{proof}

%%%%%%%%%%%%%%%%%%%%%%%%%%%%%%%%%%%%%%%%%%%%%%%%%%%%%%%%%%%%%%%%%%%%%%%%%

\end{document}